\documentclass[11pt]{article}

\usepackage[T1]{fontenc}
\usepackage[utf8]{inputenc}
\usepackage[margin=1.15in]{geometry}
\usepackage{amsmath,amssymb,amsthm}
\usepackage{graphicx}
\usepackage[authoryear,round]{natbib}
\usepackage[colorlinks=true,linkcolor=blue,citecolor=blue,urlcolor=blue]{hyperref}

\numberwithin{equation}{section}
\newtheoremstyle{imsplain}{8pt}{8pt}{\itshape}{\parindent}{\scshape}{.}{.5em}%
  {\thmname{#1}\thmnumber{ #2}\thmnote{ {\normalfont(#3)}}}
\newtheoremstyle{imsdef}{8pt}{8pt}{\normalfont}{\parindent}{\scshape}{.}{.5em}%
  {\thmname{#1}\thmnumber{ #2}\thmnote{ {\normalfont(\textit{#3})}}}
\theoremstyle{imsplain}
\newtheorem{theorem}{Theorem}[section]
\newtheorem{lemma}{Lemma}[section]
\theoremstyle{imsdef}
\newtheorem{condition}{Condition}

\newtheorem*{conditionPJ}{Condition PJ}
\newtheorem{remark}{Remark}[section]
\newtheorem{example}{Example}[section]

\newcommand{\R}{\mathbb{R}}
\newcommand{\NN}{\mathbb{N}}
\newcommand{\X}{\mathbf{X}}            % compact support set of the regressors
\newcommand{\cX}{\mathcal{X}}          % regressor draws attached to Poisson points
\newcommand{\cH}{\mathcal{H}}
\newcommand{\cZ}{\mathcal{Z}}
\newcommand{\cD}{\mathcal{D}}
\newcommand{\cI}{\mathcal{I}}
\newcommand{\cL}{\mathcal{L}}
\newcommand{\cQ}{\mathcal{Q}}
\newcommand{\cT}{\mathcal{T}}
\newcommand{\cE}{\mathcal{E}}
\newcommand{\bc}{\mathbf{c}}
\newcommand{\be}{\mathbf{e}}
\newcommand{\bN}{\mathbf{N}}           % limit Poisson process
\newcommand{\hbN}{\widehat{\mathbf{N}}}% empirical point process
\newcommand{\hZ}{\widehat{Z}}
\newcommand{\ind}{\mathbf{1}}          % indicator; use \mathds{1} (dsfont) if preferred
\newcommand{\dto}{\overset{d}{\to}}
\newcommand{\pto}{\overset{p}{\to}}
\newcommand{\argmin}{\operatorname*{arg\,min}}

\DeclareMathOperator{\Var}{Var}
\DeclareMathOperator{\Cov}{Cov}

\begin{document}

\title{EXTREMAL QUANTILE REGRESSION\thanks{Supported in part by an Alfred P.\ Sloan
Foundation Dissertation Fellowship.}}
\author{Victor Chernozhukov\\ Massachusetts Institute of Technology}
\date{}
\maketitle
\begingroup
\renewcommand{\thefootnote}{}%
\footnotetext{Originally published in \textit{The Annals of Statistics}, 2005, Vol.~33, No.~2, 806--839 (DOI: 10.1214/009053604000001165). This revised arXiv version reflects a proof audit and minor corrections conducted via Gemini Neon (November 2026). All main results (Theorems~3.1, 4.1, 5.1, and~6.1) remain unchanged; corrections address minor sign, notation, and cross-reference typos in Remarks~3.2 and~4.5, Theorem~5.1, and Section~9 (Lemmas~9.1, 9.4, 9.6, and~9.7 and Sections~9.3, 9.5, and~9.7), and Remark~4.4 adds a general-position condition ensuring a.s.\ uniqueness for Lebesgue-almost every $k>0$.}%
\addtocounter{footnote}{-1}%
\endgroup

\begin{abstract}
Quantile regression is an important tool for estimation of conditional quantiles of a
response $Y$ given a vector of covariates $X$. It can be used to measure the effect of
covariates not only in the center of a distribution, but also in the upper and lower
tails. This paper develops a theory of quantile regression in the tails. Specifically,
it obtains the large sample properties of extremal (extreme order and intermediate
order) quantile regression estimators for the linear quantile regression model with the
tails restricted to the domain of minimum attraction and closed under tail equivalence
across regressor values. This modeling setup combines restrictions of extreme value
theory with leading homoscedastic and heteroscedastic linear specifications of
regression analysis. In large samples, extreme order regression quantiles converge
weakly to $\arg\min$ functionals of stochastic integrals of Poisson processes that
depend on regressors, while intermediate regression quantiles and their functionals
converge to normal vectors with variance matrices dependent on the tail parameters and
the regressor design.
\end{abstract}


% ======================================================================
\section{Introduction}\label{sec:1}

Regression quantiles [\citet{KB78}] estimate conditional quantiles of a response
variable $Y$ given regressors $X$. They extend Laplace's \citeyearpar{Laplace} median
regression (least absolute deviation estimator) and generalize the ordinary sample
quantiles to the regression setting. Regression quantiles are used widely in empirical
work and studied extensively in theoretical statistics. See, for example,
\citet{Buchinsky}, \citet{Chamberlain}, \citet{CDS}, \citet{GJ}, \citet{HK},
\citet{Knight98}, \citet{KP87}, \citet{PK97}, \citet{Portnoy91a} and \citet{Powell},
among others.

Many potentially important applications of regression quantiles involve the study of
various extremal phenomena. In econometrics, motivating examples include the analysis of
factors that contribute to extremely low infant birthweights [cf.\ \citet{Abrevaya}];
the analysis of the highest bids in auctions [cf.\ \citet{DP}]; and estimation of
factors of high risk in finance [cf.\ \citet{Tsay} and \citet{CU}, among others]. In
biostatistics and other areas, motivating examples include the analysis of survival at
extreme durations [cf.\ \citet{KG}]; the analysis of factors that impact the approximate
boundaries of biological processes [cf.\ \citet{Cade}]; image reconstruction and other
problems where conditional quantiles near maximum or minimum are of interest [cf.\
\citet{KST}].

An important peril to inference in the listed examples is that conventional large sample
theory for quantile regression does not apply sufficiently far in the tails. In the
nonregression case, this problem is familiar, well documented and successfully dealt
with by modern extreme value theory; see, for example, \citet{LLR}, \citet{Resnick} and
\citet{EKM}. The purpose of this paper is to develop an asymptotic theory for quantile
regression in the tails based on this theory. Specifically, this paper obtains the large
sample properties of extremal (extreme order and intermediate order) quantile regression
for the class of linear quantile regression models with conditional tails of the
response variable restricted to the domain of minimum attraction and closed under the
tail equivalence across conditioning values.

The paper is organized as follows. After an introductory Section~\ref{sec:2},
Section~\ref{sec:3} joins together the linear quantile regression model with the tail
restrictions of modern extreme value theory. These restrictions are imposed in a manner
that allows regressors to impact the conditional tail quantiles of response $Y$
differently than the central quantiles. The resulting modeling setup thus covers
conventional location shift regression models, as well as more general quantile
regression models. Section~\ref{sec:4} provides the asymptotic theory for the sample
regression quantiles under the extreme order condition, $\tau_T T\to k>0$, where
$\tau_T$ is the quantile index and $T$ is the sample size. By analogy with the extreme
order quantiles in nonregression cases, the extreme order regression quantiles converge
to extreme type variates (functionals of multivariate Poisson processes that depend on
regressors). Our analysis of the case $\tau_T T\to k>0$ builds on and complements the
analysis of $\tau_T T\to 0$ given by \citet{FR}, \citet{Smith}, \citet{PJ} and
\citet{Knight01} for various types of location shift models. [\citet{C98} also studied
some nonparametric cases.] Section~\ref{sec:5} derives the asymptotic distributions of
regression quantiles under the intermediate order condition: $\tau_T T\to\infty$,
$\tau_T\to 0$, thus providing a quantile regression analog of the results on the
intermediate univariate quantiles by \citet{DdH}. As with the intermediate quantiles in
nonregression cases, the intermediate order regression quantiles, and their functionals
such as Pickands type estimators of the extreme value index, analyzed in
Section~\ref{sec:6}, are asymptotically normal with variance determined by both the tail
parameters and the regressor design. Section~\ref{sec:7} provides an illustration,
Section~\ref{sec:8} concludes, and Section~\ref{sec:9} collects the proofs.

% ======================================================================
\section{The setting}\label{sec:2}

Suppose $Y$ is the response variable in $\R$, and $X=(1,X_{-1}')'$ is a $d\times1$
vector of regressors (typically transformations of original regressors). (Throughout the
paper, given a vector $x$, $x_{-1}$ denotes $x$ without its first component $x_1$.)
Denote the conditional distribution of $Y$ given $X=x$ by $F_Y(\cdot|x)$. The present
focus is on $F_Y^{-1}(\tau|x)=\inf\{y\colon F_Y(y|x)>\tau\}$, where $\tau$ is close to
$0$. Let there be a sample
\[
\{Y_t,X_t,\ t=1,\dots,T\}\qquad\text{where }X_t\in\X,
\]
generated by a probability model with a conditional quantile function of the classical
linear-in-parameter form
\begin{equation}\label{e:2.1}
F_Y^{-1}(\tau|x)=x'\beta(\tau)\qquad\text{for all }\tau\in\cI,\ x\in\X,
\end{equation}
where $\beta(\cdot)$ is a nonparametric function of $\tau$, which when $\cI=(0,1)$ also
corresponds to the stochastic model with random coefficients:
\begin{equation}\label{e:2.2}
Y=X'\beta(\varepsilon),\qquad \varepsilon\overset{d}{=}U(0,1),\ X\in\X.
\end{equation}
Here it is necessary that \eqref{e:2.1} holds for
\begin{equation}\label{e:2.3}
\cI=[0,\eta]\qquad\text{for some }0<\eta\le1\text{ and }x\in\X,\text{ a compact subset
of }\R^d.
\end{equation}
Different linear models \eqref{e:2.1} can be applied to different covariate regions
$\X$ [which can be local neighborhoods of a given $x_0$, in which case the linear model
\eqref{e:2.1} is motivated as a Taylor expansion]. The model \eqref{e:2.1} plays a
fundamental role in the theoretical and practical literature on quantile regression
mentioned in the Introduction. Its appealing feature is the ability to capture
quantile-specific covariate effects in a convenient linear framework.

In the sequel, we combine the linear model \eqref{e:2.1} with the tail restrictions from
extreme value theory to develop applicable asymptotic results. It is of vital
consequence to impose these restrictions in a manner that preserves the quantile-specific
covariate effects, as motivated by the empirical examples listed in the Introduction.
For instance, in the analysis of U.S.\ birthweights, \citet{Abrevaya} finds that smoking
and the absence of prenatal care impact the low conditional quantiles of birthweights
much more negatively than the central birthweight quantiles. The linear framework
\eqref{e:2.1} is able to accommodate this type of impact through the quantile-specific
coefficients $\beta(\tau)$, where $\beta_{-1}(\tau)$, for $\tau$ near $0$, describes the
effect of covariate factors on extremely low birthweights and, say, $\beta_{-1}(1/2)$
describes the effect on central birthweights. Thus, when imposing extreme value
restrictions, it is important to preserve this ability.

The inference about $\beta(\tau)$ is based on the regression quantile statistics
$\hat\beta(\tau)$ [\citet{KB78}] defined by the least asymmetric absolute deviation
problem:
\begin{equation}\label{e:2.4}
\hat\beta(\tau)\in\argmin_{\beta\in\R^d}\sum_{t=1}^T\rho_\tau(Y_t-X_t'\beta)
\qquad\text{where }\rho_\tau(u)=(\tau-\ind(u\le0))u,
\end{equation}
of which Laplace's \citeyearpar{Laplace} median regression is an important case with
$\rho_{1/2}(u)=|u|/2$. The statistics $\hat\beta(\tau)$ naturally generalize the
ordinary sample quantiles to the conditional setting. In fact, the usual univariate
$\tau$-quantiles can be recovered as the solution to this problem without covariates,
that is, when $X_t=1$. [E.g., if $\tau T\in(0,1)$, $\hat\beta(\tau)=Y_{(1)}$, and if
$\tau T\in(1,2)$, $\hat\beta(\tau)=Y_{(2)}$, etc.]

In order to provide large sample properties of $\hat\beta(\tau)$ in the tails, we
distinguish three types of sample regression quantiles, following the classical theory of
order statistics: (i) an extreme order sequence, when $\tau_T\searrow0$,
$\tau_TT\to k>0$, (ii) an intermediate order sequence, when $\tau_T\searrow0$,
$\tau_TT\to\infty$, (iii) a central order sequence, when $\tau\in(0,1)$ is fixed, and
$T\to\infty$ (under which the conventional theory applies). We consider
$\hat\beta(\tau_T)$ under the extreme and intermediate order sequences, and refer to
$\hat\beta(\tau_T)$ under both sequences as the \textit{extremal regression quantiles}.
In what follows, we omit the $T$ in $\tau_T$ whenever it does not cause confusion.

% ======================================================================
\section{The extreme value restrictions on the linear quantile regression model}\label{sec:3}

This section joins the linear model \eqref{e:2.1} together with the tail restrictions
from extreme value theory, examines the consequences and presents examples.

Consider a random variable $u$ with distribution function $F_u$ and lower end-point
$s_u=0$ or $s_u=-\infty$. Recall [cf.\ \citet{Resnick}] that $F_u$ is said to have tail
of type 1, 2 or 3 if for
\begin{equation}\label{e:3.1}
\begin{array}{@{}l@{\quad}l@{\qquad}r@{}}
\text{type 1:} & \text{as } z\searrow s_u=0\text{ or }-\infty, & \\[2pt]
& F_u(z+va(z))\sim F_u(z)e^{v} & \forall\, v\in\R,\ \xi\equiv0,\\[6pt]
\text{type 2:} & \text{as } z\searrow s_u=-\infty, & \\[2pt]
& F_u(vz)\sim v^{-1/\xi}F_u(z) & \forall\, v>0,\ \xi>0,\\[6pt]
\text{type 3:} & \text{as } z\searrow s_u=0, & \\[2pt]
& F_u(vz)\sim v^{-1/\xi}F_u(z) & \forall\, v>0,\ \xi<0,
\end{array}
\end{equation}
where $a(z)\equiv\int_{s_u}^zF_u(v)\,dv/F_u(z)$, for $z>s_u$. The number $\xi$ is
commonly called the extreme value index, and $F_u$ with tails of types 1--3 is said to
belong to the domain of minimum attraction. [$a(z)\sim b(z)$ denotes that $a(z)/b(z)\to1$
as a specified limit over $z$ is taken.]

\begin{condition}\label{R1}
In addition to \eqref{e:2.1}, there exists an auxiliary line $x\mapsto x'\beta_r$ such
that for
\begin{equation}\label{e:3.2}
U\equiv Y-X'\beta_r\qquad\text{with }s_U=0\text{ or }s_U=-\infty,
\end{equation}
and some $F_u$ with type 1, 2 or 3 tails,
\begin{equation}\label{e:3.3}
F_U(z|x)\sim K(x)\cdot F_u(z)\qquad\text{uniformly in }x\in\X,\ \text{as }z\searrow s_U,
\end{equation}
where $K(\cdot)>0$ is a continuous bounded function on $\X$. Without loss of generality,
let $K(x)=1$ at $x=\mu_X$ and $F_u(z)\equiv F_U(z|\mu_X)$.
\end{condition}

\begin{condition}\label{R2}
The distribution function of $X=(1,X_{-1}')'$, $F_X$, has compact support $\X$ with
$EXX'$ positive definite. Without loss of generality, let
$\mu_X=EX=(1,0,\dots,0)'$.
\end{condition}

When $Y$ has a finite lower endpoint, that is, $X'\beta(0)>-\infty$, it is implicit in
Condition~R1 that $\beta_r\equiv\beta(0)$ so that $U\equiv Y-X'\beta(0)\ge0$ has endpoint
$0$ by construction. In the unbounded support case, $X'\beta(0)=-\infty$ and is not
suitable as an auxiliary line, but existence of any other line such that Condition~R1
holds suffices.

Condition~R1 is the main assumption. First, Condition~R1 requires the tails of
$U=Y-X'\beta_r$ for some $\beta_r$ to be in the domain of minimum attraction, which is a
nonparametric class of distributions [cf.\ \citet{Resnick} and \citet{EKM}]. In this
sense, the specification Condition~R1 is semiparametric. Examples~\ref{ex:3.1} and
\ref{ex:3.2} present some of the regression models covered by Condition~R1. Second,
Condition~R1 also requires that, for any $x',x''\in\X$, $z\mapsto F_U(z|x')$ and
$z\mapsto F_U(z|x'')$ are tail equivalent up to a constant. This condition is motivated
by the closure of the domain of minimum attraction under tail equivalence [cf.\
Proposition 1.19 in \citet{Resnick}].

Compactness of $\X$ in Condition~R1 is necessary, as the limit theory for regression
quantiles may generally change otherwise. In applications, compactness may be imposed by
the explicit trimming of observations depending on whether $X_t\in\X$. In this case the
linear model \eqref{e:2.1} is assumed to apply only to values of $X$ in $\X$. Clearly,
the smaller $\X$, the less restrictive is the linear model by virtue of Taylor
approximation [e.g., \citet{Chaudhuri91}]. Also, trimming $X$ to $\X$ eliminates the
impact of outlying values on the limit distribution and inference, as it does in the case
of the central regression quantiles. In some cases it should be possible to make $X$
unbounded by imposing higher level nonprimitive conditions, for example, similar to those
on page 98 in \citet{Knight01}. However, since we view $\X$ as a ``small'' neighborhood
over which the linear approximation \eqref{e:2.1} is adequate, we do not pursue this
extension.

Theorem~\ref{thm:3.1} shows that the function $K(x)$ in Condition~R1 can be represented
by the following types. Other properties of the linear quantile regression model under
Conditions~R1 and R2 are obtained in Lemma~\ref{lem:9.1} given in
Section~\ref{sec:9.1}.

\begin{theorem}[Three types of $K(\cdot)$]\label{thm:3.1}
Under Conditions~R1 and R2, for some $\bc\in\R^d$,
\begin{equation}\label{e:3.4}
K(x)=\begin{cases}
e^{-x'\bc}, & \text{when }F_u\text{ has type 1 tails, }\xi=0,\\
(x'\bc)^{1/\xi}, & \text{when }F_u\text{ has type 2 tails, }\xi>0,\\
(x'\bc)^{1/\xi}, & \text{when }F_u\text{ has type 3 tails, }\xi<0,
\end{cases}
\end{equation}
where $\mu_X'\bc=1$ for type 2 and 3 tails, $\mu_X'\bc=0$ for type 1 tails, and $x'\bc>0$
for all $x\in\X$ for types 2 and 3.
\end{theorem}

\begin{remark}\label{rem:3.1}
The condition $X'\bc>0$ a.s.\ for tails of types 2 and 3 arises from the linearity
assumption \eqref{e:2.1}. Indeed, \eqref{e:2.1} imposes that the quantiles should not
cross: if $l>1$, then $X'(\beta(l\tau)-\beta(\tau))>0$ a.s. Since by
Lemma~\ref{lem:9.1}(v)
$X'(\beta(l\tau)-\beta(\tau))/\mu_X'(\beta(l\tau)-\beta(\tau))\to X'\bc$ as
$\tau\searrow0$, the noncrossing condition requires $X'\bc>0$ a.s. In location-scale
shift models (cf.\ Example~\ref{ex:3.2}), the condition $X'\bc>0$ a.s.\ is equivalent to
a logical restriction on the scale function ($X'\sigma>0$ a.s.). In location shift models
(cf.\ Example~\ref{ex:3.1}), this condition is ordinarily satisfied since $X'\bc=1$ a.s.\
for tails of types 2 and 3.
\end{remark}

\begin{remark}\label{rem:3.2}
The general case when $P\{K(X)\neq1\}>0$ will be referred to as the heterogeneous case,
and $\bc$ will be referred to as the \textit{heterogeneity index}. The special case with
\begin{equation}\label{e:3.5}
K(X)=1\qquad\text{a.s.}
\end{equation}
will be referred to as the \textit{homogeneous} case. The latter amounts to $\bc=\mathbf 0$
for type 1 tails, and $\bc=\be_1\equiv(1,0,\dots)'$ for type 2 and 3 tails.
Notice that in this case $X'\bc=1$ a.s.\ for types 2 and 3 and $X'\bc=0$ a.s.\ for type 1
tails.
\end{remark}

In developing regularity conditions which target regression applications, it is natural
to try to cover the most conventional regression settings and, hopefully, more general
stochastic specifications. The following examples clarify this possibility.

\begin{example}[Location shift regression]\label{ex:3.1}
Consider the location-shift model
\begin{equation}\label{e:3.6}
Y=X'\beta+U,
\end{equation}
where $U$ is independent of $X$, and suppose $U$ is in the domain of minimum attraction.
When the lower endpoint of the support of $U$ is finite, it is normalized to $0$.
Clearly, this is a special case of Condition~R1 where $X'\beta_r\equiv X'\beta$,
$U\equiv Y-X'\beta$, $K(X)=1$ a.s. The data generating process \eqref{e:3.6} has been
widely adopted in regression work at least since \citet{Huber} and \citet{Rao}. A variety
of standard survival and duration models also imply \eqref{e:3.6} after a transformation,
for example, the Cox models with Weibull hazards and accelerated failure time models
[cf.\ \citet{DG}]. Also, \eqref{e:3.6} underlies many theoretical studies of quantile
regression. Hence, it is useful that Condition~R1 covers \eqref{e:3.6}.
\end{example}

\begin{example}[Location-scale shift regression]\label{ex:3.2}
As a generalization of \eqref{e:3.6}, consider the stochastic equation
\begin{equation}\label{e:3.7}
Y=X'\beta+X'\sigma\cdot V,\qquad V\text{ is independent of }X,
\end{equation}
where $X'\sigma>0$ (a.s.) is the scale function, and $V$ is in the domain of minimum
attraction with $\xi\neq0$. \eqref{e:3.7} implies the following linear conditional
quantile function
\begin{equation}\label{e:3.8}
F_Y^{-1}(\tau|X)=X'\beta+X'\sigma\cdot F_V^{-1}(\tau).
\end{equation}
Then for $X'\beta_r\equiv X'\beta$, $U\equiv Y-X'\beta_r=X'\sigma\cdot V$, we have
$P(X'\sigma\cdot V\le z|X)\sim(X'\sigma)^{1/\xi}\cdot F_V(z)$ as $z\searrow0$ or
$-\infty$, so Condition~R1 is satisfied with $F_u\equiv F_V$ and $K(X)=(X'\sigma)^{1/\xi}$.
The data generating process \eqref{e:3.7} has been adopted in, for example,
\citet{KB82}, \citet{GJ} and \citet{He}.
\end{example}

\begin{example}[Quantile-shift regression]\label{ex:3.3}
To see that Condition~R1 covers more general stochastic models than \eqref{e:3.6} and
\eqref{e:3.7}, note that Condition~R1 requires that $F_U(u|X)$ or $F_V(u|X)$ be
independent of $X$ only in the tails. In both cases, these weaker independence
requirements allow $X$, for example, to have a negative impact on the high and low
quantiles but to have a positive impact on the median quantiles. In contrast, notice from
\eqref{e:3.8} that \eqref{e:3.6} and \eqref{e:3.7} preclude such quantile-specific
impacts. Thus, Condition~R1 preserves the heterogeneous impact property of
\eqref{e:2.1}, allowing the impact of covariate factors on extreme quantiles to be very
different from their impact on the central quantiles.
\end{example}

% ======================================================================
\section{Asymptotics of extreme order regression quantiles}\label{sec:4}

Consider sequences $\tau_i$, $i=1,\dots,l$, such that $\tau_iT\to k_i>0$ as
$T\to\infty$, and the corresponding normalized regression quantile statistics
$\hZ_T(k_i)$, where
\begin{equation}\label{e:4.1}
\hZ_T(k)\equiv a_T(\hat\beta(\tau)-\beta_r-b_T\be_1),
\end{equation}
$\hat\beta(\tau)$ is the regression quantile, $\beta_r$ is the coefficient of the
auxiliary line defined in \eqref{e:3.2}, $\be_1\equiv(1,0,\dots)'\in\R^d$, and
$(a_T,b_T)$ are the canonical normalization constants, given by
\begin{equation}\label{e:4.2}
\begin{array}{@{}l@{\quad}l@{\qquad}l@{}}
\text{for type 1 tails:} & a_T=1/a\bigl[F_u^{-1}\bigl(\tfrac1T\bigr)\bigr], &
b_T=F_u^{-1}\bigl(\tfrac1T\bigr),\\[6pt]
\text{for type 2 tails:} & a_T=-1/F_u^{-1}\bigl(\tfrac1T\bigr), & b_T=0,\\[6pt]
\text{for type 3 tails:} & a_T=1/F_u^{-1}\bigl(\tfrac1T\bigr), & b_T=0,
\end{array}
\end{equation}
where $F_u$ is defined in Condition~R1. Moreover, consider the centered statistic
\begin{equation}\label{e:4.3}
\hZ_T^c(k)\equiv a_T(\hat\beta(\tau)-\beta(\tau))
\end{equation}
and the point process, for $U_t=Y_t-X_t'\beta_r$,
\begin{equation}\label{e:4.4}
\hbN(\cdot)=\sum_{t=1}^T\ind(\{a_T(U_t-b_T),X_t\}\in\cdot).
\end{equation}
We will show that $\hbN(\cdot)$ converges weakly to the Poisson process
\begin{equation}\label{e:4.5}
\bN(\cdot)=\sum_{i=1}^\infty\ind(\{J_i,\cX_i\}\in\cdot),
\end{equation}
with points $\{J_i,\cX_i\}$ satisfying
\begin{equation}\label{e:4.6}
(J_i,\cX_i,\ i\ge1)=\begin{cases}
(\ln(\Gamma_i)+\cX_i'\bc,\ \cX_i), & \text{for type 1 tails,}\\
(-\Gamma_i^{-\xi}\cX_i'\bc,\ \cX_i), & \text{for type 2 tails,}\\
(\Gamma_i^{-\xi}\cX_i'\bc,\ \cX_i), & \text{for type 3 tails,}
\end{cases}\qquad i\ge1,
\end{equation}
where $\{\cX_i\}$ is an i.i.d.\ sequence with law $F_X$,
\begin{equation}\label{e:4.7}
\Gamma_i\equiv\sum_{j=1}^i\cE_j,\qquad i\ge1,
\end{equation}
and $\{\cE_j\}$ is an i.i.d.\ sequence of unit-exponential variables, independent of
$\{\cX_i\}$. In the homogeneous case \eqref{e:3.5}, $J_i$ and $\cX_i$ are independent
since
\begin{equation}\label{e:4.8}
\cX_i'\bc=\begin{cases}0, & \text{for type 1 tails,}\\ 1, & \text{for type 2 and 3 tails,}
\end{cases}\qquad\text{for all }i\ge1.
\end{equation}
The following theorem establishes the weak limit of $\hZ_T(k)$'s as a function of $\bN$.

\begin{theorem}[Extreme order regression quantiles]\label{thm:4.1}
Assume Conditions~R1 and R2 and that $\{Y_t,X_t\}$ is an i.i.d.\ or a stationary sequence
satisfying the Meyer type conditions of Lemma~\ref{lem:9.4}. Then as $\tau_T\to k>0$ and
$T\to\infty$,
\begin{equation}\label{e:4.9}
\hZ_T(k)\dto Z_\infty(k)\equiv\argmin_{z\in\cZ}\Bigl[-k\mu_X'z+\int(x'z-u)^+\,d\bN(u,x)\Bigr],
\end{equation}
provided $Z_\infty(k)$ is a uniquely defined random vector in $\cZ$, where
$(x'z-u)^+=\ind(u\le x'z)(x'z-u)$, $\cZ=\R^d$ for type 1 and 3 tails, and
$\cZ=\{z\in\R^d\colon\max_{x\in\X}z'x\le0\}$ for type 2 tails. Moreover,
\begin{equation}\label{e:4.10}
\hZ_T^c(k)\dto Z_\infty^c(k)\equiv Z_\infty(k)-\eta(k),
\end{equation}
where
\begin{equation}\label{e:4.11}
\eta(k)=\begin{cases}
\bc+\ln k\,\be_1, & \text{for type 1 tails,}\\
-k^{-\xi}\bc, & \text{for type 2 tails,}\\
k^{-\xi}\bc, & \text{for type 3 tails.}
\end{cases}
\end{equation}
If $Z_\infty(k)$ is a uniquely defined random vector for $k=k_1,\dots,k_l$,
\begin{align*}
(\hZ_T(k_1)',\dots,\hZ_T(k_l)')'&\dto(Z_\infty(k_1)',\dots,Z_\infty(k_l)')',\\
(\hZ_T^c(k_1)',\dots,\hZ_T^c(k_l)')'&\dto(Z_\infty^c(k_1)',\dots,Z_\infty^c(k_l)')'.
\end{align*}
\end{theorem}

\begin{remark}[The limit criterion function]\label{rem:4.1}
The limit objective function $-k\mu_X'z+\int(x'z-u)^+\,d\bN(u,x)$ can also be written as
\begin{equation}\label{e:4.12}
-k\mu_X'z+\sum_{i=1}^\infty(\cX_i'z-J_i)^+.
\end{equation}
\end{remark}

\begin{remark}[Homogeneous case]\label{rem:4.2}
The limit result is simpler for the homogeneous case \eqref{e:3.5}, since $\bN$ does not
depend on the heterogeneity parameter $\bc$ due to \eqref{e:4.8}.
\end{remark}

\begin{remark}[Case with $\tau T\to0$]\label{rem:4.3}
The linear programming estimator, which corresponds to $T\tau\to0$ in \eqref{e:2.4} (in
comparison, here $\tau T\to k>0$), was studied in \citet{FR}, \citet{Smith}, \citet{PJ},
\citet{Knight99,Knight01} and \citet{C98} under various types of location-shift
specification \eqref{e:3.6}. This estimator is the solution to the problem
\begin{equation}\label{e:4.13}
\max_{\beta\in\R^d}\bar X'\beta\ \text{ such that }\ Y_t\ge X_t'\beta
\quad\text{for all }t\le T,\qquad \bar X=T^{-1}\sum_{t=1}^TX_t.
\end{equation}
The asymptotics of \eqref{e:4.13} and proofs differ substantively from the ones given here
for $\tau T\to k>0$. The analysis of $\tau T\to k>0$ is specifically motivated by the
applications listed in the Introduction.
\end{remark}

\begin{remark}[Uniqueness]\label{rem:4.4}
The limit objective function is convex, and it is assumed in Theorem~\ref{thm:4.1} that
$Z_\infty(k)$ is unique and tight. Lemma~\ref{lem:9.7} shows that a sufficient condition
for tightness is the design condition of \citet{PJ}. Taking tightness as given,
conditions for uniqueness can be established. Define $\cH$ as the set of all $d$-element
subsets of $\NN$. For $h\in\cH$, let $\cX(h)$ and $J(h)$ be the matrix with rows
$\cX_t$, $t\in h$, and vector with elements $J_t$, $t\in h$, respectively. Let
$\cH^*=\{h\in\cH\colon|\cX(h)|\ne0\}$. $\cH^*$ is nonempty a.s.\ by Condition~R2 and is
countable. Application of the argument of Theorem 3.1 of \citet{KB78} gives that an
$\arg\min$ of \eqref{e:4.12} takes the form $z_h=\cX(h)^{-1}J(h)$ for some $h\in\cH^*$,
and must satisfy the gradient condition
\begin{equation}\label{e:4.14}
\zeta_k(z_h)\equiv\Bigl(k\mu_X-\sum_{t=1}^\infty\ind(J_t<\cX_t'z_h)\cX_t\Bigr)'\cX(h)^{-1}
\in[0,1]^d,
\end{equation}
where the $\arg\min$ is unique iff $\zeta_k(z_h)\in\cD=(0,1)^d$. Thus, uniqueness holds
for a fixed $k>0$ if
\begin{equation}\label{e:4.15}
P(\zeta_k(z_h)\in\partial\cD\text{ for some }h\in\cH^*)=0.
\end{equation}
This condition is a direct analog of Koenker and Bassett's \citeyearpar{KB78} condition
for uniqueness in finite samples; for instance, it is satisfied for a given $k$ when
covariates $\cX_{-1t}$ are absolutely continuous [cf.\ \citet{Portnoy91b}]. Thus,
uniqueness holds generically in the sense that for a fixed $k$ adding arbitrarily small
absolutely continuous perturbations to $\{\cX_{-1t}\}$ ensures \eqref{e:4.15}.

Moreover, without requiring absolute continuity of $\cX_{-1t}$, \eqref{e:4.15} holds for
Lebesgue-almost every $k>0$ under the general-position condition
$P(\mu_X\in\operatorname{span}(\cX_1,\dots,\cX_{d-1}))=0$ (for $d\ge2$; for $d=1$,
$k\neq\lceil k\rceil$ suffices as in Remark~\ref{rem:4.6}). Indeed, for each $h\in\cH^*$
and finite subset $\mu\subset\NN\setminus h$, the $r$-th coordinate of
$(k\mu_X-\sum_{t\in\mu}\cX_t)'\cX(h)^{-1}$ is affine in $k$ with slope $\mu_X'v_r$, where
$v_r$ is the $r$-th column of $\cX(h)^{-1}$. Since $\mu_X'v_r=0$ iff $\mu_X$ lies in the
span of the remaining $d-1$ rows of $\cX(h)$, the condition implies $\mu_X'v_r\ne0$ a.s.\
for all $h\in\cH^*$ and $r$, so each of the countably many triples $(h,\mu,r)$ hits
$\{0,1\}$ at most twice as $k$ varies; Fubini's theorem then yields \eqref{e:4.15} for
Lebesgue-almost every $k>0$. For $d=2$ with $X=(1,X_{-1})'$, this condition reduces simply
to $P(X_{-1}=E[X_{-1}])=0$, which holds for any binary regressor with $P(X_{-1}=1)\in(0,1)$
(whereas a symmetric three-point regressor $X_{-1}\in\{-1,0,1\}$ fails it); for $d\ge3$, it
means that $E[X_{-1}]$ does not lie on an affine hyperplane in $\R^{d-1}$ spanned by $d-1$
support points with positive probability.
\end{remark}

\begin{remark}[Asymptotic density]\label{rem:4.5}
The density of $Z_\infty(k)$ can be stated following \citet{KB78}. Given $\{\cX_t\}$,
$h\in\cH^*$, and $J(h)$, the probability that $Z_\infty(k)=\cX(h)^{-1}J(h)$ equals
$P\{\zeta_k(\cX(h)^{-1}J(h))\in\cD\,|\,\{\cX_t\},J(h)\}$. Conditional on
$\{Z_\infty(k)=\cX(h)^{-1}J(h)\}$, $h\in\cH^*$, and $\cX(h)$, the density of
$Z_\infty(k)$ at $z$ is $f_{J(h)|\cX(h)}(\cX(h)z)\cdot|\cX(h)|$, where
$f_{J(h)|\cX(h)}(u)$, $u\in\R^d$, is the joint density of $J(h)$ conditional on
$\cX(h)$. Thus, the joint density of $Z_\infty(k)$ at $z$ is
\[
f_{Z_\infty(k)}(z)=E\Bigl[\sum_{h\in\cH^*}f_{J(h)|\cX(h)}(\cX(h)z)\cdot|\cX(h)|
\times P\{\zeta_k(z)\in\cD\,|\,\{\cX_t\},J(h)=\cX(h)z\}\Bigr].
\]
Finally, for $f_{Z_\infty(k)}(z)$ to be nondefective, $Z_\infty(k)=O_p(1)$ should be
established (cf.\ Lemma~\ref{lem:9.7}).
\end{remark}

\begin{remark}[Univariate case]\label{rem:4.6}
The density simplifies in the classical nonregression case, that is, when $X=1$, in which
case we also have the simplification \eqref{e:4.8}. In this case, an $\arg\min$ is
necessarily an order statistic, that is, $z_h=J(h)=J_h$; the gradient condition
\eqref{e:4.14} becomes
\begin{equation}\label{e:4.16}
\zeta_k(z_h)\equiv\Bigl(k-\sum_{t=1}^\infty\ind(J_t<z_h)\Bigr)\in[0,1];
\end{equation}
and the condition for uniqueness is that $\zeta_k(z_h)\in\cD=(0,1)$. Then, for
$k\ne\lceil k\rceil$, $P\{\zeta_k(z_h)\in\cD\}=1$ if $h=\lceil k\rceil$ and
$P\{\zeta_k(z_h)\in\cD\}=0$ if $h\ne\lceil k\rceil$. Here $k\ne\lceil k\rceil$ is needed
for uniqueness. Hence, $f_{Z_\infty(k)}(z)=f_{J_{\lceil k\rceil}}(z)$, which is the limit
density of the $\lceil k\rceil$th order statistics in the univariate case. Thus,
uniqueness holds for almost every $k\in(0,\infty)$.
\end{remark}

% ======================================================================
\section{Asymptotics of intermediate order regression quantiles}\label{sec:5}

In order to develop asymptotic results for the intermediate regression quantiles, the
following additional Condition~R3 will be added. First, existence of the quantile density
function $\partial F_U^{-1}(\tau|x)/\partial\tau\equiv x'\partial\beta(\tau)/\partial\tau$
and its regular variation will be required. Second, the tail equivalence of the
conditional distribution functions, previously assumed in Condition~R1, will now be
strengthened to the tail equivalence of conditional quantile density functions.

\begin{condition}\label{R3}
In addition to Conditions~R1 and R2, for $\xi$ defined in \eqref{e:3.1},
\begin{equation}\label{e:5.1}
\begin{array}{@{}l@{\ \ }l@{}}
\text{(i)} & \dfrac{\partial F_U^{-1}(\tau|x)}{\partial\tau}\sim
\dfrac{\partial F_u^{-1}(\tau/K(x))}{\partial\tau}\qquad\text{uniformly in }x\in\X,\\[12pt]
\text{(ii)} & \dfrac{\partial F_u^{-1}(\tau)}{\partial\tau}\ \text{ is regularly varying at
$0$ with exponent }-\xi-1.
\end{array}
\end{equation}
\end{condition}

In the homoscedastic case \eqref{e:3.5}, Condition~R3(i) amounts to
$\frac{\partial F_U^{-1}(\tau|x)}{\partial\tau}\sim\frac{\partial F_u^{-1}(\tau)}{\partial\tau}$
uniformly in $x\in\X$ as $\tau\searrow0$. Condition~R3(ii) is a von Mises type condition;
see \citet{DdH} for a detailed analysis of the plausibility of Condition~R3(ii).

For an intermediate sequence such that $\tau\searrow0$ and $\tau T\to\infty$, define, for
$m>1$,
\begin{equation}\label{e:5.2}
\hZ_T\equiv a_T(\hat\beta(\tau)-\beta(\tau)),\qquad
a_T\equiv\frac{\sqrt{\tau T}}{\mu_X'(\beta(m\tau)-\beta(\tau))}.
\end{equation}
Consider also $k$ sequences $\{\tau l_1,\dots,\tau l_k\}$, where $l_1,\dots,l_k$ are
positive constants, and corresponding statistics $(\hZ_T(l_1)',\dots,\hZ_T(l_k)')'$,
where, for $l>0$ and $m>1$,
\begin{equation}\label{e:5.3}
\hZ_T(l)\equiv a_T(l)(\hat\beta(l\tau)-\beta(l\tau)),\qquad
a_T(l)\equiv\frac{\sqrt{\tau lT}}{\mu_X'(\beta(ml\tau)-\beta(l\tau))}.
\end{equation}
The following theorem establishes the weak limits for $\hZ_T$ and $\hZ_T(l)$'s. Because
$\tau\searrow0$, the limits depend only on the tail parameters $\xi$ and $\bc$, as in
Theorem~\ref{thm:4.1}, but since $\tau T\to\infty$, the limits are normal, unlike in
Theorem~\ref{thm:4.1}.

\begin{theorem}[Intermediate order regression quantiles]\label{thm:5.1}
Suppose Conditions~R1--R3 hold, and that $\{Y_t,X_t\}$ is an i.i.d.\ sequence or a
stationary series satisfying the conditions of Lemma~\ref{lem:9.6}. Then, as
$\tau T\to\infty$ and $\tau\searrow0$,
\begin{equation}\label{e:5.4}
\hZ_T\dto Z_\infty=N(0,\Omega_0),\qquad
\Omega_0\equiv\cQ_H^{-1}\cQ_X\cQ_H^{-1}\frac{\xi^2}{(m^{-\xi}-1)^2},
\end{equation}
where, for $\xi=0$, interpret $\xi^2/(m^{-\xi}-1)^2$ as $(\ln m)^{-2}$ and
\begin{align}
\cQ_H&\equiv E[H(X)^{-1}XX'],\qquad \cQ_X\equiv EXX',\label{e:5.5}\\
H(x)&\equiv x'\bc\ \text{ for type 2 and 3 tails},\qquad H(x)\equiv1\ \text{ for type 1 tails}.\label{e:5.6}
\end{align}
In addition,
\begin{align}
(\hZ_T(l_1)',\dots,\hZ_T(l_k)')'&\dto(Z_\infty(l_1)',\dots,Z_\infty(l_k)')'=N(0,\Omega),\label{e:5.7}\\
EZ_\infty(l_i)Z_\infty(l_j)'&=\Omega_0\times\min(l_i,l_j)/\sqrt{l_il_j}.\label{e:5.8}
\end{align}
Finally, $a_T(l)$ can be replaced by $\sqrt{\tau lT}/\bar X'(\hat\beta(ml\tau)-\hat\beta(l\tau))$
without affecting \eqref{e:5.4} and \eqref{e:5.7}, that is,
\begin{equation}\label{e:5.9}
a_T(l)\Big/\biggl(\frac{\sqrt{\tau lT}}{\bar X'(\hat\beta(ml\tau)-\hat\beta(l\tau))}\biggr)
\pto1\qquad\text{where }\bar X=T^{-1}\sum_{t=1}^TX_t.
\end{equation}
\end{theorem}

\begin{remark}[Scaling constants]\label{rem:5.1}
It may be useful to have the same normalization $a_T$ in place of $a_T(l)$ for the joint
convergence. This is possible by noting that $a_T/a_T(l)\to l^{-\xi}/\sqrt l$.
\end{remark}

\begin{remark}[Homogeneous case]\label{rem:5.2}
In the homogeneous case \eqref{e:3.5}, $H(X)=1$, so the variance simplifies to
\begin{equation}\label{e:5.10}
\Omega_0=\cQ_X^{-1}\frac{\xi^2}{(m^{-\xi}-1)^2}.
\end{equation}
\end{remark}

\begin{remark}[Nonregression case]\label{rem:5.3}
Theorem~\ref{thm:5.1} extends Theorem 3.1 of \citet{DdH}, which applies to univariate
quantiles, to the case of regression quantiles. In fact, Theorem 3.1 of \citet{DdH} can be
specialized from Theorem~\ref{thm:5.1} with $X=1$ and $m=2$. In this case the variance
becomes
\begin{equation}\label{e:5.11}
\frac{\xi^2}{(2^{-\xi}-1)^2}=\frac{2^{2\xi}\xi^2}{(2^{\xi}-1)^2},
\end{equation}
as \citet{DdH} found in their Theorem 3.1.
\end{remark}

% ======================================================================
\section{Quantile regression spacings and tail inference}\label{sec:6}

The tail parameters enter the limit distributions in Theorems~\ref{thm:4.1} and
\ref{thm:5.1}, and estimation of the tail index is an important problem of its own. The
following results show how to estimate them by applying \citet{Pickands} type procedures
to the quantile regression spacings.

Consider the following parameters and statistics:
\begin{equation}\label{e:6.1}
\begin{aligned}
\varphi&=\frac{x'(\hat\beta(m\tau)-\hat\beta(\tau))}{x'(\beta(m\tau)-\beta(\tau))},\\[4pt]
\rho_{x,\dot x,l}&=\frac{x'(\beta(ml\tau)-\beta(l\tau))}{\dot x'(\beta(m\tau)-\beta(\tau))},\\[4pt]
\hat\rho_{x,\dot x,l}&=\frac{x'(\hat\beta(ml\tau)-\hat\beta(l\tau))}{\dot x'(\hat\beta(m\tau)-\hat\beta(\tau))}.
\end{aligned}
\end{equation}
Theorem~\ref{thm:6.1} shows that the quantile regression spacings of intermediate order
consistently approximate the corresponding spacings in the population [results (i) and
(ii)], which then reveal the tail parameters [results (iii) and (iv)].

\begin{theorem}[Quantile regression spacings and tail inference]\label{thm:6.1}
Suppose the conditions of Theorem~\ref{thm:5.1} hold. Then as $\tau\searrow0$,
$\tau T\to\infty$, for all $l>0$, $m>1$, $x,\dot x\in\X$,
\begin{enumerate}
\item[(i)] $\varphi\pto1$,
\item[(ii)] $\hat\rho_{x,\dot x,l}-\rho_{x,\dot x,l}\pto0$,\quad
$\rho_{x,\dot x,l}\to l^{-\xi}\cdot[H(x)/H(\dot x)]$, for $H(x)$ defined in
Theorem~\ref{thm:5.1},
\item[(iii)] $\hat\xi_{rp}\equiv\dfrac{-1}{\ln l}\ln\hat\rho_{\bar X,\bar X,l}\pto\xi$,
\item[(iv)] $\hat\rho_{x,\bar X,1}\pto x'\bc$ uniformly in $x\in\X$ ($\xi\ne0$),
\item[(v)] for $\pi=\mu_X'\cQ_H^{-1}\cQ_X\cQ_H^{-1}\mu_X$, $l=m=2$, if
$\sqrt{\tau T}(\rho_{\bar X,\bar X,l}-\lim_T\rho_{\bar X,\bar X,l})\to0$,
\begin{equation}\label{e:6.2}
\sqrt{\tau T}(\hat\xi_{rp}-\xi)\dto N\Bigl(0,\ \pi\cdot\frac{\xi^2(2^{2\xi+1}+1)}{(2(2^\xi-1)\ln2)^2}\Bigr).
\end{equation}
\end{enumerate}
\end{theorem}

\begin{remark}[Homogeneous case]\label{rem:6.1}
The proposed estimator $\hat\xi_{rp}$ consistently estimates the tail index $\xi$ in the
heteroscedastic and homoscedastic quantile regression models, and it is a regression
extension of the \citet{Pickands} estimator. In fact, in the homoscedastic model
\eqref{e:3.5} or when $X=1$, $\pi=\mu_X'(EXX')^{-1}\mu_X=\be_1'(EXX')^{-1}\be_1=1$, so the
variance in \eqref{e:6.2} reduces to that of the canonical Pickands estimator.
\end{remark}

% ======================================================================
\section{An illustrative example}\label{sec:7}

The set of results established here may provide reliable and practical inference for
extremal regression quantiles. To illustrate this possibility, the following simple
example compares graphically the conventional central asymptotic approximation, where,
for fixed $\tau\in(0,1)$ as $T\to\infty$,
\begin{equation}\label{e:7.1}
\sqrt T(\hat\beta(\tau)-\beta(\tau))\dto N\Bigl(0,\ \frac{1}{f_U^2(F_U^{-1}(\tau))}(EXX')^{-1}\tau(1-\tau)\Bigr),
\end{equation}
to the extreme approximation (cf.\ Theorem~\ref{thm:4.1}). The comparison is based on the
following design: $\tau=0.025$, $Y_t=X_t'\beta+U_t$, $U_t\sim\text{Cauchy}$,
$t=1,\dots,500$, where $X_t=(1,X_{-1t}')'\in\R^5$, $X_{-1t}$ are i.i.d.\ Beta$(3,3)$
variables, and $\beta=(1,1,1,1,1)$. [A more detailed simulation study is given in
\citet{C99}.] In this comparison, the parameters of the limit distribution are fixed at
the true values.

\begin{figure}[t]
\centering
\includegraphics[width=\textwidth]{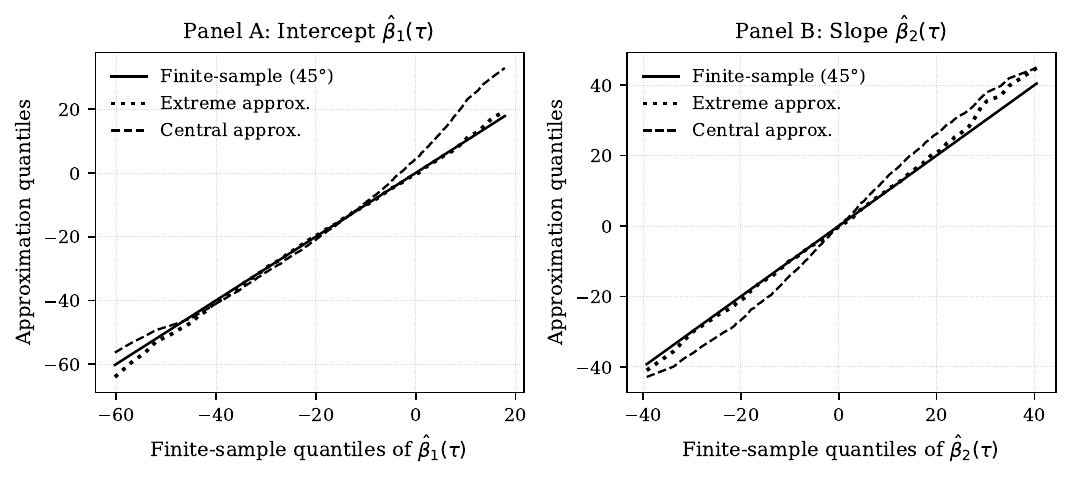}
\caption{Panel A plots quantiles of the finite-sample distribution of
$\hat\beta_1(\tau)$ (horizontal axis) against the quantiles of the extreme approximation
(cf.\ Theorem~\ref{thm:4.1}) and the quantiles of the central approximation
\eqref{e:7.1} (vertical axis). Panel B plots quantiles of the finite-sample distribution
of $\hat\beta_2(\tau)$ (horizontal axis) against the quantiles of the extreme
approximation (cf.\ Theorem~\ref{thm:4.1}) and the quantiles of the central approximation
\eqref{e:7.1} (vertical axis). The plot is based on 3,000 simulations of the regression
model described in Section~\ref{sec:7} (regenerated in vector format for this revised version). The dashed line ``- - - -'' denotes quantiles of
the central approximation, and the dotted line ``$\cdot\cdot\cdot\cdot\cdot\cdot$''
denotes quantiles of the extreme approximation (this approximation almost coincides with
``------''). The simulated quantiles of the finite-sample distribution are given by the
45-degree line depicted as the solid line ``------.''}
\label{fig:1}
\end{figure}

Figure~\ref{fig:1} plots (a) quantiles of the simulated finite-sample distribution of
$\hat\beta_1(0.025)$ and $\hat\beta_2(0.025)$, (b) quantiles of the simulated extreme
approximation (cf.\ Theorem~\ref{thm:4.1}), (c) quantiles of the central approximation
[cf.\ \eqref{e:7.1}]. Here $\tau\times T=0.025\times500=12.5$. It can be seen that the
extreme approximation accurately captures the actual sampling distribution of both the
intercept estimator $\hat\beta_1(0.025)$ and the slope estimator $\hat\beta_2(0.025)$. In
contrast, the central approximation \eqref{e:7.1} does not capture asymmetry and thick
tails of the true finite sample distribution. The intermediate approximation (cf.\
Theorem~\ref{thm:5.1}), performs similarly to the central approximation and is not
plotted. The central and intermediate approximations are expected to perform better for
less extreme quantiles.

% ======================================================================
\section{Conclusion}\label{sec:8}

The paper obtains the large sample properties of extreme order and intermediate order
quantile regression for the class of linear quantile regression models with tails of the
response variable restricted to the domain of minimum attraction and closed under tail
equivalence across conditioning values. There are several interesting directions for
future work. It would be important to determine the most practical and reliable inference
procedures that can be based on the obtained limit distributions. Also, it would be
interesting to examine estimation of the extreme conditional quantiles defined through an
extrapolation of the intermediate regression quantiles. The nonregression case has been
considered in \citet{DdH} and \citet{dHR}, and the approach may prove useful in the
quantile regression case. Another interesting direction would be an investigation of the
Hill and other tail index estimators based on regression quantiles.

% ======================================================================
\section{Proofs}\label{sec:9}

% ----------------------------------------------------------------------
\subsection{Properties of the linear quantile regression model under Conditions R1 and R2}\label{sec:9.1}
Let
\begin{align}
M&\equiv\text{any fixed compact sub-interval of }(0,1)\cup(1,\infty),\label{e:9.1}\\
M'&\equiv\text{any other fixed compact sub-interval of }(0,1)\cup(1,\infty),\label{e:9.2}\\
\cT(\tau')&\equiv\{\tau\colon\tau=s\tau',\ s\in\cL\}\qquad\text{where }\tau'\searrow0,\label{e:9.3}\\
\cL&\equiv\text{any fixed compact sub-interval of }(0,\infty).\label{e:9.4}
\end{align}

\begin{lemma}[Properties of the linear model under Conditions R1 and R2]\label{lem:9.1}
Conditions R1 and R2 imply that \textup{(}for a constant vector $\bc$ specified in
Theorem~\textup{\ref{thm:3.1})}:
\begin{enumerate}
\item[(i)] $K(x)$ can be represented by the forms specified in Theorem~\ref{thm:3.1}.
\item[(ii)] $a_T(\beta(\tau)-\beta_r-b_T\be_1)\to\eta(k)$ for $\eta(k)$ defined in
Theorem~\ref{thm:4.1}.
\item[(iii)] Uniformly in $(m,\tau,x)\in M\times\cT(\tau')\times\X$, as $\tau'\searrow0$,
\begin{equation}\label{e:9.5}
\frac{\beta_{-1}(\tau)-\beta_{-1r}}{F_u^{-1}(m\tau)-F_u^{-1}(\tau)}\to\mu(m)=
\begin{cases}
\dfrac{\bc_{-1}}{m^{-\xi}-1}, & \text{for }\xi<0,\\[8pt]
\dfrac{\bc_{-1}}{m^{-\xi}-1}, & \text{for }\xi>0,\\[8pt]
\dfrac{\bc_{-1}}{\ln m}, & \text{for }\xi=0;
\end{cases}
\end{equation}
also $\beta_1(\tau)-\beta_{1r}=F_u^{-1}(\tau)$, and
$(\beta_{-1}(\tau)-\beta_{-1r})/F_u^{-1}(\tau)\to\bc_{-1}$ for $\xi\ne0$.
\item[(iv)] Uniformly in $(m,\tau,x)\in M\times\cT(\tau')\times\X$, as $\tau'\searrow0$,
\begin{equation}\label{e:9.6}
\frac{(x-\mu_X)'(\beta(\tau)-\beta_r)}{\mu_X'(\beta(m\tau)-\beta(\tau))}\to
\begin{cases}
(x-\mu_X)'\dfrac{\bc}{m^{-\xi}-1}, & \text{if }\xi<0,\\[8pt]
(x-\mu_X)'\dfrac{\bc}{m^{-\xi}-1}, & \text{if }\xi>0,\\[8pt]
(x-\mu_X)'\dfrac{\bc}{\ln m}, & \text{if }\xi=0.
\end{cases}
\end{equation}
\item[(v)] Uniformly in $(m,\tau,x)\in M\times\cT(\tau')\times\X$, as $\tau'\searrow0$,
\begin{equation}\label{e:9.7}
\frac{x'(\beta(m\tau)-\beta(\tau))}{\mu_X'(\beta(m\tau)-\beta(\tau))}\to
\begin{cases}
x'\bc, & \text{if }\xi<0,\\
x'\bc, & \text{if }\xi>0,\\
1, & \text{if }\xi=0.
\end{cases}
\end{equation}
\item[(vi)] Uniformly in $(l,m,\tau,x)\in M\times M'\times\cT(\tau')\times\X$, as
$\tau'\searrow0$,
\begin{equation}\label{e:9.8}
\frac{x'(\beta(l\tau)-\beta(\tau))}{x'(\beta(m\tau)-\beta(\tau))}\to
\begin{cases}
\dfrac{l^{-\xi}-1}{m^{-\xi}-1}, & \text{if }\xi<0,\\[8pt]
\dfrac{1-l^{-\xi}}{1-m^{-\xi}}, & \text{if }\xi>0,\\[8pt]
\dfrac{\ln l}{\ln m}, & \text{if }\xi=0.
\end{cases}
\end{equation}
\end{enumerate}
\end{lemma}

Write $F_u\in D(H_\xi)$ if $F_u$ is a c.d.f.\ in the domain of minimum attraction with tail
index $\xi$. Write $F_u\in\mathcal R_\gamma(0)$ if $F_u$ is a regularly varying function at
$0$ with exponent $\gamma$.

\begin{lemma}[Useful relations]\label{lem:9.2}
Under Conditions R1 and R2, uniformly in $(m,l,\tau)\in M\times M'\times\cT(\tau')$, as
$\tau'\searrow0$:
\begin{enumerate}
\item[(i)] Suppose $F_1(z)\sim F_2(z)$ as $z\searrow0$ or $-\infty$ and $F_1\in D(H_\xi)$.
Then $F_2\in D(H_\xi)$; $F_1^{-1}$ and $F_2^{-1}\in\mathcal R_{-\xi}(0)$;
$F_1(F_1^{-1}(\tau))\sim\tau$ and $F_2(F_2^{-1}(\tau))\sim\tau$; and
\begin{equation}\label{e:9.9}
(F_1^{-1}(m\tau)-F_1^{-1}(\tau))\sim(F_2^{-1}(m\tau)-F_2^{-1}(\tau)).
\end{equation}
\item[(ii)] If $F_U(z|x)\sim K(x)F_u(z)$ as $z\searrow0$ or $-\infty$ for each $x\in\X$
\textup{(}compact\textup{)}, where $K(x)\in(0,\infty)$ for all $x\in\X$, then for each
$x\in\X$,
\begin{equation}\label{e:9.10}
F_U^{-1}(m\tau|x)-F_U^{-1}(\tau|x)\sim F_u^{-1}(m\tau/K(x))-F_u^{-1}(\tau/K(x)).
\end{equation}
\item[(iii)] $\displaystyle\frac{F_u^{-1}(m\tau)-F_u^{-1}(\tau)}{F_u^{-1}(l\tau)-F_u^{-1}(\tau)}
\to\begin{cases}\dfrac{m^{-\xi}-1}{l^{-\xi}-1} & \text{if }\xi<0,\\[8pt]
\dfrac{1-m^{-\xi}}{1-l^{-\xi}} & \text{if }\xi>0,\\[8pt]
\dfrac{\ln m}{\ln l} & \text{if }\xi=0;\end{cases}$\quad for $F_u\in D(H_\xi)$.
\item[(iv)] $\displaystyle\frac{F_u^{-1}(lm\tau)-F_u^{-1}(l\tau)}{a(F_u^{-1}(\tau))}\to\ln m$
if $F_u\in D(H_0)$, where $a(\cdot)$ is the auxiliary function defined in \eqref{e:3.1}.
\end{enumerate}
\end{lemma}

\begin{proof}[Proof]
Results (i), (iii) and (iv) are well known [cf.\ \citet{dH84} and \citet{Resnick}, Chapters
1 and 2]. Result (ii) holds from (i) pointwise in $x$.
\end{proof}

\begin{proof}[Proof of Lemma~\ref{lem:9.1}]
Claim (i): The proof consists of two steps, where we use notation
$(\cL,M,\cT(\tau'),\tau')$ as defined in \eqref{e:9.1}\allowbreak--\eqref{e:9.4}.

\medskip\noindent\textsc{Step 1.} In this step all of the results hold uniformly in
$(m,\tau,x)\in M\times\cT(\tau')\times\X$ as $\tau'\searrow0$, but we shall suppress this
qualification for notational simplicity. By construction in Condition~R1,
$x'(\beta(\tau)-\beta_r)\equiv F_U^{-1}(\tau|x)$ and
$\mu_X'(\beta(\tau)-\beta_r)\equiv F_U^{-1}(\tau|\mu_X)\equiv F_u^{-1}(\tau)$. Hence,
\begin{equation}\label{e:9.11}
B_\tau(x,m)\equiv\frac{(x-\mu_X)'(\beta(\tau)-\beta_r)}{\mu_X'(\beta(m\tau)-\beta(\tau))}
\equiv\frac{F_U^{-1}(\tau|x)-F_u^{-1}(\tau)}{F_u^{-1}(m\tau)-F_u^{-1}(\tau)}.
\end{equation}
We would like to show that, for each $x\in\X$,
\begin{equation}\label{e:9.12}
B_\tau(x,m)\to B(x,m)\equiv
\begin{cases}
\dfrac{(1/K(x))^{-\xi}-1}{m^{-\xi}-1}, & \text{if }\xi<0,\\[8pt]
\dfrac{1-(1/K(x))^{-\xi}}{1-m^{-\xi}}, & \text{if }\xi>0,\\[8pt]
\dfrac{\ln(1/K(x))}{\ln m}, & \text{if }\xi=0.
\end{cases}
\end{equation}
We will show \eqref{e:9.12} for the case $\xi<0$ only; others follow similarly. Fix any
$x\in\X$. By Condition~R1 and Lemma~\ref{lem:9.2}(i), $F_U(F_U^{-1}(\tau|x)|x)\sim\tau$.
Hence, by Condition~R1, $K(x)\cdot F_u(F_U^{-1}(\tau|x))\sim\tau$ as $\tau'\searrow0$.
Therefore, there exist sequences of constants $K_\tau(x)$ and $K_\tau'(x)$ such that
\begin{equation}\label{e:9.13}
\begin{gathered}
F_u^{-1}(\tau/K_\tau(x))\le F_U^{-1}(\tau|x)\le F_u^{-1}(\tau/K_\tau'(x))\\
\text{where }K_\tau(x)\to K(x)\text{ and }K_\tau'(x)\to K(x).
\end{gathered}
\end{equation}
Therefore,
\begin{equation}\label{e:9.14}
\begin{aligned}
\frac{F_u^{-1}(\tau/K_\tau(x))-F_u^{-1}(\tau)}{F_u^{-1}(m\tau)-F_u^{-1}(\tau)}
&\le B_\tau(x,m)\\
&\le\frac{F_u^{-1}(\tau/K_\tau'(x))-F_u^{-1}(\tau)}{F_u^{-1}(m\tau)-F_u^{-1}(\tau)}.
\end{aligned}
\end{equation}
Suppose that $K(x)\ne1$. By Lemma~\ref{lem:9.2}(iii),
\begin{equation}\label{e:9.15}
\frac{F_u^{-1}(\tau/K_\tau(x))-F_u^{-1}(\tau)}{F_u^{-1}(m\tau)-F_u^{-1}(\tau)}
\to\frac{(1/K(x))^{-\xi}-1}{m^{-\xi}-1}=B(x,m),
\end{equation}
and, likewise, conclude for $K_\tau'(x)$ in place of $K_\tau(x)$. Therefore,
$B_\tau(x,m)\to B(x,m)$ when $K(x)\ne1$. To show that $B_\tau(x,m)\to B(x,m)$ also holds
for $K(x)=1$ with $B(x,m)=0$, let $\kappa'$ and $\kappa''$ be any positive constants such
that $\kappa'<1<\kappa''$. By monotonicity of the quantile function, for all sufficiently
small $\tau'$,
\begin{equation}\label{e:9.16}
\begin{aligned}
\frac{F_u^{-1}(\tau/\kappa'')-F_u^{-1}(\tau)}{F_u^{-1}(m\tau)-F_u^{-1}(\tau)}
&\le\frac{F_u^{-1}(\tau/K_\tau(x))-F_u^{-1}(\tau)}{F_u^{-1}(m\tau)-F_u^{-1}(\tau)}\\
&\le\frac{F_u^{-1}(\tau/\kappa')-F_u^{-1}(\tau)}{F_u^{-1}(m\tau)-F_u^{-1}(\tau)}.
\end{aligned}
\end{equation}
By Lemma~\ref{lem:9.2}(iii), as $\tau'\searrow0$, the upper and lower bounds in \eqref{e:9.16}
converge to
\begin{equation}\label{e:9.17}
\frac{(1/\kappa'')^{-\xi}-1}{m^{-\xi}-1}\qquad\text{and}\qquad
\frac{(1/\kappa')^{-\xi}-1}{m^{-\xi}-1}.
\end{equation}
If in \eqref{e:9.17} we let $\kappa',\kappa''\to1$, then expressions in \eqref{e:9.17}
$\to0$. Therefore, since $\kappa'$ and $\kappa''$ can be chosen arbitrarily close to $1$,
it follows from \eqref{e:9.16} and \eqref{e:9.17} that
$\frac{F_u^{-1}(\tau/K_\tau(x))-F_u^{-1}(\tau)}{F_u^{-1}(m\tau)-F_u^{-1}(\tau)}\to0$ as
$\tau'\searrow0$. Likewise, conclude for $K_\tau'(x)$ in place of $K_\tau(x)$. Therefore,
$B_\tau(x,m)\to B(x,m)=0$ when $K(x)=1$.

\medskip\noindent\textsc{Step 2.} By Step 1, for each $x\in\X$, uniformly in
$(m,\tau)\in M\times\cT(\tau')$ as $\tau'\searrow0$,
\begin{equation}\label{e:9.18}
B_\tau(x,m)=\frac{(x-\mu_X)'(\beta(\tau)-\beta_r)}{\mu_X'(\beta(m\tau)-\beta(\tau))}\to B(x,m).
\end{equation}
Since (a) $B(x,m)$ is finite and continuous in $x$ over $\X$ by conditions imposed on
$K(x)$ in Condition~R1, and (b) $B_\tau(x,m)$ is linear in $x$, the relation
\eqref{e:9.18} also holds uniformly in $x\in\X$. Recall that $(x-\mu_X)_1=0$. Since
$(x-\mu_X)_{-1}$ ranges over a nondegenerate subset of $\R^{d-1}$, \eqref{e:9.18} implies
\begin{equation}\label{e:9.19}
\frac{\beta_{-1}(\tau)-\beta_{-1r}}{\mu_X'(\beta(m\tau)-\beta(\tau))}\to\mu(m),
\end{equation}
uniformly in $(m,\tau)\in M\times\cT(\tau')$ as $\tau'\searrow0$, where $\mu(m)$ is some
vector of finite constants. Hence, $B(x,m)$ is affine in $(x-\mu_X)$. Note also that
$(x-\mu_X)=(0,x_{-1}')'$. Therefore, if $\xi=0$, $B(x,m)$ affine and
$B(x,m)=-\ln K(x)/\ln m$ imply $K(x)=e^{-(x-\mu_X)'\bc}=e^{-x_{-1}'\bc_{-1}}=e^{-x'\bc}$ for
all $x$ iff $\bc_1=0$.
When $\xi<0$, $B(x,m)$ affine and $B(x,m)=(K(x)^\xi-1)/(m^{-\xi}-1)$ imply
$K(x)=(1+(x-\mu_X)'\bc)^{1/\xi}$, which equals $(x'\bc)^{1/\xi}$ for all $x$ iff
$\bc_1=1$. Likewise, conclude for $\xi>0$. This completes the proof of claim (i).

Claim (iii) follows directly from \eqref{e:9.19} and the preceding paragraph.

Claim (iv) is verified by substituting the forms of $K(x)$ found above into
\eqref{e:9.18}.

Claim (v) holds pointwise in $x$ by Lemma~\ref{lem:9.2}(ii) and (iii). Since the
left-hand side in \eqref{e:9.7} is linear in $x$ and $\X$ is compact, it also holds
uniformly in $x\in\X$.

A combination of Lemma~\ref{lem:9.2}(iii) with claim (v) implies claim (vi).

Claim (ii). If $\xi<0$, by claim (iii) uniformly in $k$ in any compact subset of
$(0,\infty)$ as $T\to\infty$,
\begin{equation}\label{e:9.20}
\begin{aligned}
&a_T\Bigl(\beta\Bigl(\frac kT\Bigr)-\beta_r\Bigr)\\
&\qquad\sim a_T\bc F_u^{-1}\Bigl(\frac kT\Bigr)=\bc F_u^{-1}\Bigl(\frac kT\Bigr)\Big/F_u^{-1}\Bigl(\frac1T\Bigr)\to k^{-\xi}\bc,
\end{aligned}
\end{equation}
since by Lemma~\ref{lem:9.2}(i) $F_u^{-1}\in\mathcal R_{-\xi}(0)$; similarly, if $\xi>0$,
\begin{equation}\label{e:9.21}
\begin{aligned}
&a_T\Bigl(\beta\Bigl(\frac kT\Bigr)-\beta_r\Bigr)\\
&\qquad\sim a_T\bc F_u^{-1}\Bigl(\frac kT\Bigr)=-\bc F_u^{-1}\Bigl(\frac kT\Bigr)\Big/F_u^{-1}\Bigl(\frac1T\Bigr)\to -k^{-\xi}\bc.
\end{aligned}
\end{equation}
If $\xi=0$, by $\bc_1=0$, Lemma~\ref{lem:9.2}(i), (iv) and claim (iii) [using $m=e$ in
$\mu(m)$], we have that uniformly in $k$ in any compact subset of $(0,\infty)$,
\begin{equation}\label{e:9.22}
\begin{aligned}
&a_T\Bigl(\beta\Bigl(\frac kT\Bigr)-\beta_r-b_T\be_1\Bigr)\\
&\quad\sim\frac{1}{a(F^{-1}(1/T))}\\
&\qquad\times\Bigl[\bc\Bigl(F_u^{-1}\Bigl(e\frac kT\Bigr)-F_u^{-1}\Bigl(\frac kT\Bigr)\Bigr)
+\be_1\Bigl(F_u^{-1}\Bigl(\frac kT\Bigr)-F_u^{-1}\Bigl(\frac1T\Bigr)\Bigr)\Bigr]\\
&\quad\to\bc\ln e+\be_1\ln k=\bc+\be_1\ln k.
\end{aligned}
\end{equation}
\end{proof}

% ----------------------------------------------------------------------
\subsection{Proof of Theorem~\ref{thm:3.1}}\label{sec:9.2}
Follows from Lemma~\ref{lem:9.1}(i).\hfill$\square$

% ----------------------------------------------------------------------
\subsection{Proof of Theorem~\ref{thm:4.1}}\label{sec:9.3}

\medskip\noindent\textit{Part 1.} Referring to \eqref{e:2.4}, notice that $Z_T(k)$ defined
in \eqref{e:4.1} solves
\begin{equation}\label{e:9.23}
\hZ_T(k)\in\argmin_{z\in\R^d}\Bigl[\frac1{a_T}\sum_{t=1}^T\rho_\tau(a_T(U_t-b_T)-X_t'z)\Bigr]
\end{equation}
[where $z\equiv a_T(\beta-\beta_r-b_T\be_1)$]. Rearranging terms, the objective function
becomes
\begin{equation}\label{e:9.24}
\begin{aligned}
\frac1{a_T}\Bigl[&-\tau T\bar X'z-\sum_{t=1}^T\ind(a_T(U_t-b_T)\le X_t'z)(a_T(U_t-b_T)-X_t'z)\\
&+\tau\cdot\sum_{t=1}^T a_T(U_t-b_T)\Bigr].
\end{aligned}
\end{equation}
Multiply \eqref{e:9.24} by $a_T$ and subtract
\begin{equation}\label{e:9.25}
\sum_{t=1}^T\ind(a_T(U_t-b_T)\le-\delta)(-\delta-a_T(U_t-b_T))+\sum_{t=1}^T\tau a_T(U_t-b_T)
\qquad\text{for some }\delta>0,
\end{equation}
which does not affect optimization, and denote the new objective function $Q_T(z,k)$:
\begin{equation}\label{e:9.26}
Q_T(z,k)\equiv-\tau T\bar X'z+\sum_{t=1}^Tl_\delta(a_T(U_t-b_T),X_t'z),
\end{equation}
where
\begin{equation}\label{e:9.27}
l_\delta(u,v)\equiv\ind(u\le v)(v-u)-\ind(u\le-\delta)(-\delta-u)\qquad\text{for }\delta>0.
\end{equation}
Since it is a sum of convex functions in $z$, $Q_T(z,k)$ is convex in $z$. The
transformations make (as shown later) $Q_T$ a continuous functional of the point process
$\hbN$:
\begin{equation}\label{e:9.28}
Q_T(z,k)=-\tau T\bar X'z+\int_El_\delta(j,x'z)\,d\hbN(j,x),
\end{equation}
where the point process
\begin{equation}\label{e:9.29}
\hbN(\cdot)\equiv\sum_{t\le T}\ind\{(a_T(U_t-b_T),X_t)\in\cdot\}
\end{equation}
is taken to be a random element of the metric space $M_p(E)$ of point processes defined on
the measure space $(E,\mathcal E)$ and equipped with the metric induced by the topology of
vague convergence [cf.\ \citet{Resnick}].

It will suffice to restrict our attention to underlying measure spaces $(E,\mathcal E)$ of
the form
\begin{equation}\label{e:9.30}
E=\begin{cases}
E_1\equiv[-\infty,\infty)\times\X, & \text{for type 1 tails,}\\
E_2\equiv[-\infty,0)\times\X, & \text{for type 2 tails,}\\
E_3\equiv[0,\infty)\times\X, & \text{for type 3 tails,}
\end{cases}
\end{equation}
with $\sigma$-algebra $\mathcal E$ generated by the open sets of $E$. The topology on
$E_1$, $E_2$ and $E_3$ is assumed to be standard so that, for example,
$[-\infty,a]\times\X$ is compact in $E_2$ for $a<0$ and in $E_1$ for any $a<\infty$.

Part 2 shows that, for type 1 and 3 tails, the marginal weak limit of $Q_T$ is a finite
convex function in $z$:
\begin{equation}\label{e:9.31}
Q_\infty(z,k)=-k\mu_X'z+\int_El_\delta(j,x'z)\,d\bN(j,x),\qquad z\in\R^d,
\end{equation}
where $\bN$ is the Poisson point process defined in the statement of Theorem~\ref{thm:4.1}.

Part 2 also shows that, for type 2 tails, the marginal weak limit of $Q_T$ is a finite
convex function in $z$:
\begin{equation}\label{e:9.32}
\begin{gathered}
Q_\infty(z,k)=-k\mu_X'z+\int_El_\delta(j,x'z)\,d\bN(j,x)\\
\text{for }z\in\cZ_N\equiv\Bigl\{z\in\R^d\colon\max_{x\in\X}x'z<0\Bigr\},
\end{gathered}
\end{equation}
where $\bN$ is the Poisson point process defined in the statement of Theorem~\ref{thm:4.1},
and
\begin{equation}\label{e:9.33}
Q_\infty(z,k)=+\infty\qquad\text{for }z\in\cZ_P\equiv\Bigl\{z\in\R^d\colon\max_{x\in\X}x'z>0\Bigr\}.
\end{equation}
The function $Q_\infty(z,k)$ is convex and $l_\delta(j,x'z)=(x'z-j)^+\ge0$ when $j\ge-\delta$.
Hence, $Q_\infty(z,k)$ is also well defined over entire
$\cZ=\{z\in\R^d\colon\max_{x\in\X}x'z\le0\}$, although it may equal $+\infty$ at
$z\colon\max_{x\in\X}x'z=0$. Also, note that $\cZ_N\cup\cZ_P$ is dense in $\R^d$.

Recall the convexity lemma [cf.\ \citet{Geyer} and \citet{Knight99}], which states:
Suppose (i) a sequence of convex lower-semicontinous functions
$Q_T\colon\R^d\to\bar{\R}$ marginally converges to $Q_\infty\colon\R^d\to\bar{\R}$ over a
dense subset of $\R^d$, (ii) $Q_\infty$ is finite over a nonempty open set $\cZ_0$, and
(iii) $Q_\infty$ is uniquely minimized at a random vector $Z_\infty$. Then any $\arg\min$
of $Q_T$, denoted $\hZ_T$, converges in distribution to $Z_\infty$.

We showed (i) and (ii) in Part 2, and we assumed (iii). (A sufficient condition for
uniqueness is given in Remark~\ref{rem:4.4}.) Hence, application of the convexity lemma to
our case gives
\begin{equation}\label{e:9.34}
\hZ_T(k)\dto Z_\infty(k)\equiv\argmin_{z\in\R^d}Q_\infty(z,k).
\end{equation}
Note also that, for type 2, tails, the $\arg\min Z_\infty(k)$ necessarily belongs to
\[
\cZ=\Bigl\{z\in\R^d\colon\max_{x\in\X}x'z\le0\Bigr\}.
\]
This gives us the conclusion stated in
Theorem~\ref{thm:4.1} upon noting that $Q_\infty(z,k)$ differs from the limit objective
function of Theorem~\ref{thm:4.1} only by a finite random variable that does not depend on
$z$.

\medskip\noindent\textit{Part 2.} It remains to verify that (I) there exists a nonempty
open set $\cZ_0$ such that $Q_\infty(z,k)$ is finite a.s.\ for all $z\in\cZ_0$ and (II)
$Q_\infty(\cdot,k)$ is, indeed, the weak marginal limit of $Q_T(\cdot,k)$.

To show (I), when tails are of type 1 and 3, choose $\cZ_0$ as any open bounded subset of
$\R^d$; when tails are of type 2, additionally require $\cZ_0\subset\cZ_N$ for each $l$
(possible by compactness of $\X$). For any $z\in\cZ_0$, $(u,x)\mapsto l_\delta(u,x'z)$ is
in $C_K(E)$ (continuous functions on $E$ vanishing outside a compact set $K$) by the
arguments in (II). This implies $\int_El_\delta(u,x'z)\,d\bN(u,x)$ is finite a.s., since
$\bN\in M_p(E)$.

To show (II), $Q_\infty(\cdot,k)$ is the marginal weak limit of $\{Q_T(\cdot,k)\}$ iff for
any finite collection $(z_j,\ j=1,\dots,l)$,
$(Q_T(z_j,k),\ j=1,\dots,l)\dto(Q_\infty(z_j,k),\ j=1,\dots,l)$. Since
$\bar X'z_j\pto\mu_X'z_j$ and $\tau T\to k>0$, it remains to verify
\begin{equation}\label{e:9.35}
\begin{gathered}
\Bigl(\int_El_\delta(u,x'z_j)\,d\hbN(u,x),\ j=1,\dots,l\Bigr)\\
\dto\Bigl(\int_El_\delta(u,x'z_j)\,d\bN(u,x),\ j=1,\dots,l\Bigr).
\end{gathered}
\end{equation}
Define the mapping $\mathbf T\colon M_p(E)\to\R^l$ (for $E=E_1$, $E_2$ or $E_3$) by
\begin{equation}\label{e:9.36}
\mathbf T\colon N\mapsto\Bigl(\int_El_\delta(u,x'z_j)\,dN(u,x),\ j=1,\dots,l\Bigr).
\end{equation}

(a) Consider type 1 tails and set $E=E_1$. The map $(u,x)\mapsto l_\delta(u,x'z_j)$ is in
$C_K(E_1)$ (continuous functions on $E_1$ vanishing outside a compact set $K$), since by
construction it is continuous on $E_1$ and vanishes outside
$K\equiv[-\infty,\max(\kappa,-\delta)]\times\X$, where
$\kappa=\max_{x\in\X,z\in\{z_1,\dots,z_l\}}x'z$. $K$ is compact in $E_1$ since
$\kappa<\infty$ by Condition~R2. Hence, $\hbN\mapsto\mathbf T(\hbN)$ is continuous from
$M_p(E_1)$ to $\R^l$. Thus, $\hbN\Rightarrow\bN$ in $M_p(E_1)$ implies
$\mathbf T(\hbN)\dto\mathbf T(\bN)$.

(b) Consider type 3 tails and set $E=E_3$. The map $(u,x)\mapsto l_\delta(u,x'z_j)$ is in
$C_K(E_3)$: by construction, it is continuous on $E_3$ and vanishes outside
$K\equiv[0,\max(\kappa,0)]\times\X$, where $\kappa=\max_{x\in\X,z\in\{z_1,\dots,z_l\}}x'z$.
$K$ is compact in $E_3$ since $\kappa<\infty$ by Condition~R2. Therefore,
$\hbN\mapsto\mathbf T(\hbN)$ is continuous from $M_p(E_3)$ to $\R^l$. Hence,
$\hbN\Rightarrow\bN$ in $M_p(E_3)$ implies $\mathbf T(\hbN)\dto\mathbf T(\bN)$.

(c) Consider type 2 tails and set $E=E_2$. (c)(i) shows that \eqref{e:9.35} holds on
$\cZ_N$, while (c)(ii) shows that $Q_T(z,k)\pto\infty$ for any $z\in\cZ_P$. [Sets $\cZ_N$ and
$\cZ_P$ are defined in \eqref{e:9.32} and \eqref{e:9.33}.]

(i) The map $(u,x)\mapsto l_\delta(u,x'z)$ is in $C_K(E_2)$ if $z\in\cZ_N$, since, by
construction, it is continuous on $E_2$ and vanishes outside
$K\equiv[-\infty,\max(\kappa,-\delta)]\times\X$, where
$\kappa=\max_{x\in\X,z\in\{z_1,\dots,z_l\}}x'z$. $K$ is compact in $E_2$ since $\kappa<0$
if $z\in\cZ_N$. Hence, $\hbN\mapsto\mathbf T(\hbN)$ is continuous from $M_p(E_2)$ to $\R^l$.
Then $\hbN\Rightarrow\bN$ in $M_p(E_2)$ implies $\mathbf T(\hbN)\dto\mathbf T(\bN)$.

(ii) Observe that $I\equiv\sum_{t\le T}l_\delta(a_TU_t,X_t'z)\ind(a_TU_t\le-\delta)=O_p(1)$ by
the argument in (i). Observe that $l_\delta(u,v)=(v-u)^+\ge0$ for any $u\ge-\delta$. Hence,
\begin{equation}\label{e:9.37}
\begin{gathered}
l_\delta(u,v)=\ind(-\delta\le u\le v)(v-u)\ge\ind(-\delta\le u\le0,\ v\ge\varepsilon)\varepsilon\\
\text{for any }u\ge-\delta\text{ and any }\varepsilon>0.
\end{gathered}
\end{equation}
For a given $z\in\cZ_P$, since $\X$ equals the support of $X$, $\max_{x\in\X}x'z>0$ implies
that $X'z\ge\varepsilon$ occurs with positive probability for some $\varepsilon>0$. Fix
this $\varepsilon$. Since $1/a_T\to\infty$ for type 2 tails,
$P(-\delta/a_T\le U\le0,\ X'z\ge\varepsilon)\to\pi=P(U\le0,\ X'z\ge\varepsilon)>0$.
$\pi>0$ because $\inf_{x\in\X}P(U\le0|X=x)>0$ for type 2 tails by assumptions in
Condition~R1. Therefore,
$II\equiv\sum_{t\le T}\ind(-\delta/a_T\le U_t\le0,\ X_t'z\ge\varepsilon)\varepsilon\pto+\infty$
in $\bar{\R}$. Since $Q_T(z,k)\ge-k\mu_X'z+I+II$ by \eqref{e:9.37}, $Q_T(z,k)\pto+\infty$
for any $z\in\cZ_P$.

\medskip\noindent\textit{Part 3.} By Lemma~\ref{lem:9.1}(ii),
$a_T(\beta(\tau)-\beta_r-b_T\be_1)\to\eta(k)$. Hence,
$\hZ_T^c(k)\dto Z_\infty^c(k)\equiv Z_\infty(k)-\eta(k)$.

\medskip\noindent\textit{Part 4.}
$(\hZ_T(k_j)',\ j=1,\dots,l)'\in\argmin_{z\in\R^{d\times l}}[Q_T(z_1,k_1)+\cdots+Q_T(z_l,k_l)]$,
for $z=(z_1,\dots,z_l)$. Since this objective is a sum of objective functions in Parts 1
and 2, the previous derivation of the marginal limit and subsequent arguments apply very
similarly to $Q_T(z_1,k_1)+\cdots+Q_T(z_l,k_l)$ to conclude that
$(\hZ_T(k_j)',\ j=1,\dots,l)'\dto(\hZ_\infty(k_j)',\ j=1,\dots,l)'=
\argmin_{z\in\R^{d\times l}}[Q_\infty(z_1,k_1)+\cdots+Q_\infty(z_l,k_l)]$.\hfill$\square$

% ----------------------------------------------------------------------
\subsection{Weak limit of $\hbN$}\label{sec:9.4}

\begin{lemma}[{\citet[Proposition 3.22]{Resnick}}]\label{lem:9.3}
Suppose $\bN$ is a simple point process in $M_p(E)$, $\cT$ is a basis of relatively
compact open sets such that $\cT$ is closed under finite unions and intersections and, for
any $F\in\cT$, $P(\bN(\partial F)=0)=1$. Then $\hbN\Rightarrow\bN$ in $M_p(E)$ if, for all
$F\in\cT$,
\begin{align}
\lim_{T\to\infty}P[\hbN(F)=0]&=P[\bN(F)=0],\label{e:9.38}\\
\lim_{T\to\infty}E\hbN(F)&=E\bN(F)<\infty.\label{e:9.39}
\end{align}
\end{lemma}

\begin{remark}\label{rem:9.1}
In our case, $\cT$ consists of finite unions and intersections of bounded open rectangles
in $E_1$, $E_2$ and $E_3$ [cf.\ \citet{Resnick}].
\end{remark}

We impose \citet{Meyer} conditions on the ``rare'' events
$A_t^T(F)\equiv\{w\in\Omega\colon(a_T(U_t-b_T),X_t)\in F\}$.

\begin{lemma}[Poisson limits under Meyer mixing conditions]\label{lem:9.4}
Suppose that, for any $F\in\cT$, the triangular sequence of events
$\{(A_t^T(F),\ t\le T),\ T\ge1\}$ is stationary and $\alpha$-mixing with mixing coefficient
$\alpha_T(\cdot)$, condition \eqref{e:9.39} holds, and the Meyer type condition holds:
There exist sequences of integers $(p_n,\ n\ge1)$, $(q_n,\ n\ge1)$,
$(t_n=n(p_n+q_n),\ n\ge1)$ such that as $n\to\infty$, for some $r\ge1$,
(a) $n^r\alpha_{t_n}(q_n)\to0$, (b) $q_n/p_n\to0$, $p_{n+1}/p_n\to1$, and
(c) $I_{p_n}=\sum_{i=1}^{p_n-1}(p_n-i)P(A_1^{t_n}(F)\cap A_{i+1}^{t_n}(F))=o(1/n)$. Then in
$M_p(E)$, $\hbN\Rightarrow\bN$, a Poisson point process with mean measure $m\colon
m(F)\equiv\lim_{T\to\infty}E\hbN(F)$.
\end{lemma}

\begin{proof}[Proof]
For any $F\colon m(F)>0$, $\lim_{T\to\infty}P[\hbN(F)=0]=P[\bN(F)=0]=e^{-m(F)}$, by
\citet{Meyer}. The same also holds for $F\colon m(F)=0$, since $E\hbN(F)\to0$ implies
$P(\hbN(F)=0)\to1$. Conclude by Lemma~\ref{lem:9.3}.
\end{proof}

\begin{remark}\label{rem:9.2}
Condition $I_{p_n}=o(1/n)$ prevents clusters of ``rare'' events $A_t^T(F)$, eliminating
compound Poisson processes as limits.
\end{remark}

\begin{lemma}[Limit $\bN$ under Conditions R1 and R2]\label{lem:9.5}
Suppose Conditions R1 and R2 hold and that $(Y_t,X_t)$ is an i.i.d.\ or stationary
strongly mixing sequence that satisfies the conditions of Lemma~\ref{lem:9.4} with
$(a_T,b_T)$ defined in \eqref{e:4.2}. Then:
\begin{enumerate}
\item[(i)] $\hbN\Rightarrow\bN$ in $M_p(E)$, where $E=E_1,E_2$ and $E_3$ for tails of types
1, 2 and 3, respectively. $\bN$ is a Poisson point process with mean intensity measure:
$m(du,dx)=K(x)\times dh(u)\times dF_X(x)$, where $h(u)=e^u$ for type 1,
$h(u)=(-u)^{-1/\xi}$ for type 2, and $h(u)=u^{-1/\xi}$ for type 3 tails.
\item[(ii)] Points $(J_i,\cX_i)$ of $\bN$ have the representation
$(J_i,\cX_i,\ i\ge1)\overset{d}{=}(h^{-1}(\Gamma_i/K(\cX_i)),\cX_i,\ i\ge1)$, where
$h^{-1}$ is the inverse of $h$, $\Gamma_i=\cE_1+\cdots+\cE_i$, $i\ge1$ ($\{\cE_i\}$ are
i.i.d.\ standard exponential), and $\{\cX_i\}$ are i.i.d.\ r.v.s with law $F_X$,
independent of $\{\cE_i\}$.
\end{enumerate}
\end{lemma}

\begin{proof}
To show (i), by Lemmas~\ref{lem:9.3} and \ref{lem:9.4} the proof reduces to verifying
$\lim_TE\hbN(F)=m(F)$ for all $F$ in $\cT$. For example, as in Leadbetter, Lindgren and
Rootz\'en [\citeyear{LLR}, page 103], it suffices to consider $F$ of the form
$F=\bigcup_{j=1}^kF_j$, where $F_j=(l_j,u_j)\times\X_j$, where $F_1,\dots,F_k$ are
nonoverlapping, nonempty subsets of $E$, and $\X_1,\dots,\X_k$ are intersections of open
bounded rectangles of $\R^d$ with $\X$. Then by the stationarity and $F_j$'s
nonoverlapping,
\begin{align*}
E\hbN(F)&=E\sum_{t=1}^T\ind[(a_T(U_t-b_T),X_t)\in F]\\
&=\sum_{j=1}^kT\,P[(a_T(U-b_T),X)\in(l_j,u_j)\times\X_j]
\end{align*}
\begin{align}
&=\sum_{j=1}^kT\cdot E(P[(a_T(U-b_T),X)\in(l_j,u_j)\times\X_j\,|\,X])\nonumber\\
&=\sum_{j=1}^kT\cdot E(P[(a_T(U-b_T)\in(l_j,u_j)\,|\,X]\cdot\ind[X\in\X_j])\nonumber\\
&=\sum_{j=1}^kT\cdot E\bigl((F_U[u_j/a_T+b_T\,|\,X]\nonumber\\
&\hspace{6em}-F_U[l_j/a_T+b_T\,|\,X])\cdot\ind[X\in\X_j]\bigr).\label{e:9.40}
\end{align}
Suppose that $l_j>-\infty$ for all $j$. Then as $T\to\infty$,
\begin{align}
E\hbN(F)&=\sum_{j=1}^kE\Biggl(\Biggl(\frac{F_U[u_j/a_T+b_T\,|\,X]}{F_u[u_j/a_T+b_T]}\cdot T\cdot F_u[u_j/a_T+b_T]\nonumber\\
&\hspace{4em}-\frac{F_U[l_j/a_T+b_T\,|\,X]}{F_u[l_j/a_T+b_T]}\cdot T\cdot F_u[l_j/a_T+b_T]\Biggr)\cdot\ind[X\in\X_j]\Biggr)\nonumber\\
&\sim\sum_{j=1}^kE\bigl((K(X)[h(u_j)-h(l_j)])\ind[X\in\X_j]\bigr)\label{e:9.41}\\
&=\sum_{j=1}^k\int_{F_j}K(x)\,dh(u)\times dF_X(x)\nonumber\\
&=\sum_{j=1}^km(F_j)=m(F).\nonumber
\end{align}
In \eqref{e:9.41}, $\sim$ follows from two observations. First, the assumed tail
equivalence Condition~R1 implies
\begin{equation}\label{e:9.42}
\frac{F_U[l/a_T+b_T\,|\,x]}{F_u[l/a_T+b_T]}\sim K(x)\qquad\text{uniformly in }x\in\X,
\end{equation}
since by definition of $(a_T,b_T)$ given in \eqref{e:4.2}, $l/a_T+b_T\searrow F_u^{-1}(0)=0$
or $=-\infty$ for any $l\in(-\infty,\infty)$ for type 1 tails, any $l\in(-\infty,0)$ for
type 2 tails, and $l\in[0,\infty)$ for type 3 tails. Second, for example, as in
Leadbetter, Lindgren and Rootz\'en [\citeyear{LLR}, page 103], the definition of the tail
types \eqref{e:3.1} implies that (a) for tails of type 2, for any $l<0$,
$TF_u(l/a_T)=TF_u(-lF_u^{-1}(\tfrac1T))\sim(-l)^{-1/\xi}TF_u(F_u^{-1}(\tfrac1T))\sim(-l)^{-1/\xi}$,
(b) for tails of type 3, for any $l>0$,
$TF_u(l/a_T)=TF_u(lF_u^{-1}(\tfrac1T))\sim l^{-1/\xi}TF_u(F_u^{-1}(\tfrac1T))\sim l^{-1/\xi}$
and (c) for tails of type 1, for any $l\in\R$,
$TF_u(l/a_T+b_T)=TF_u(l/a(F_u^{-1}(\tfrac1T))+F_u^{-1}(\tfrac1T))\sim e^lTF_u(F_u^{-1}(\tfrac1T))\sim e^l$.

On the other hand, if for some $j$'s, $l_j=-\infty$ for type 1 or 2 tails, then we have the
replacement $TF_U[l_j/a_T+b_T|X]=0$ in \eqref{e:9.40}, and \eqref{e:9.41} follows
similarly.

To show (ii), construct a Poisson random measure (PRM) with the given $m(\cdot)$. First,
define a canonical homogeneous PRM $\bN_1$ with points $\{\Gamma_i,\ i\ge1\}$. It has the
mean measure $m_1(du)=du$ on $[0,\infty)$, for example, \citet{Resnick}. Second, by
Proposition 3.8 in \citet{Resnick}, the composed point process $\bN_2$ with points
$\{\Gamma_i,\cX_i\}$ is PRM with mean measure $m_2(du,dx)=du\times dF_X(x)$ on
$[0,\infty)\times\X$, because $\{\cX_i\}$ are i.i.d.\ and are independent of $\{\Gamma_i\}$.
Finally, the point process $\bN$ with the transformed points $\{\mathbf T(\Gamma_i,\cX_i)\}$,
where $\mathbf T\colon(u,x)\mapsto(h^{-1}(u/K(x)),x)$, is PRM with the desired mean measure
on $E\times\X$, $m(dj,dx)=m_2\circ\mathbf T^{-1}(dj,dx)=K(x)\times dh(j)\times dF_X(x)$, by
Proposition 3.7 in \citet{Resnick}.
\end{proof}

% ----------------------------------------------------------------------
\subsection{Proof of Theorem~\ref{thm:5.1}}\label{sec:9.5}
Step 1 outlines the overall proof using standard convexity arguments, while the
\textit{main} Step 2 invokes regular variation assumptions on the conditional quantile
density to demonstrate a quadratic approximation of the criterion function. Step 3 shows
joint convergence of several regression quantile statistics. Step 4 demonstrates that
$a_T$ can be estimated consistently.

\medskip\noindent\textsc{Step 1.} With reference to \eqref{e:2.4}, notice that
$\hZ_T\equiv a_T(\hat\beta(\tau)-\beta(\tau))$, defined in \eqref{e:5.2}, minimizes
\begin{equation}\label{e:9.43}
Q_T(z,\tau)\equiv\frac{a_T}{\sqrt{\tau T}}\sum_{t=1}^T\Bigl(\rho_\tau\Bigl(Y_t-X_t'\beta(\tau)-\frac{X_t'z}{a_T}\Bigr)-\rho_\tau(Y_t-X_t'\beta(\tau))\Bigr).
\end{equation}
Using Knight's identity,
\begin{equation}\label{e:9.44}
\begin{aligned}
&\rho_\tau(u-v)-\rho_\tau(u)\\
&\qquad=-v(\tau-\ind(u<0))+\int_0^v(\ind(u\le s)-\ind(u\le0))\,ds,
\end{aligned}
\end{equation}
write, a.s.,
\begin{equation}\label{e:9.45}
\begin{aligned}
Q_T(z,\tau)&=W_T(\tau)'z+G_T(z,\tau),\\
W_T(\tau)&\equiv\frac{-1}{\sqrt{\tau T}}\sum_{t=1}^T(\tau-\ind[Y_t<X_t'\beta(\tau)])X_t,\\
G_T(z,\tau)&\equiv\frac{a_T}{\sqrt{\tau T}}\Biggl(\sum_{t=1}^T\int_0^{X_t'z/a_T}\bigl[\ind(Y_t-X_t'\beta(\tau)\le s)\\
&\hspace{9em}-\ind(Y_t-X_t'\beta(\tau)\le0)\bigr]\,ds\Biggr).
\end{aligned}
\end{equation}
By Lemma~\ref{lem:9.6}, $W_T(\tau)\dto W\equiv N(0,EXX')$, and by Step 2,
\begin{equation}\label{e:9.46}
G_T(z,\tau)\pto\frac12\Bigl(\frac{m^{-\xi}-1}{-\xi}\Bigr)z'\cQ_Hz,\qquad m>1,
\end{equation}
where $\cQ_H\equiv E[H(X)^{-1}XX']$, $H(x)\equiv x'\bc$ for type 2 and 3 tails, and
$H(x)=1$ for type 1 tails. Thus, the weak marginal limit of $Q_T(z)$ is given by
\begin{equation}\label{e:9.47}
Q_\infty(z)=W'z+\frac12\cdot\Bigl(\frac{m^{-\xi}-1}{-\xi}\Bigr)\cdot z'\cQ_Hz.
\end{equation}
We have that $EXX'$ is positive definite and by Theorem~\ref{thm:3.1} that
$0<H(X)<c<\infty$ for some constant $c$. Thus, $\cQ_H$ is finite and $\cQ_H$ is positive
definite. Indeed, $z'\cQ_Hz=E(X'z)^2/H(X)=0$ for some $z\ne0$ if and only if $X'z=0$ a.s.,
which contradicts $EXX'$ positive definite. Thus, the marginal limit $Q_\infty(z)$ is
uniquely minimized at
$Z_\infty\equiv\bigl(\frac{\xi}{m^{-\xi}-1}\bigr)\cQ_H^{-1}W=
N\bigl(0,\frac{\xi^2}{(m^{-\xi}-1)^2}\cQ_H^{-1}EXX'\cQ_H^{-1}\bigr)$. By the convexity
lemma [e.g., \citet{Geyer} and \citet{Knight99}], $\hZ_T\dto Z_\infty$.

\medskip\noindent\textsc{Step 2.} This step demonstrates that as $\tau\searrow0$,
\begin{equation}\label{e:9.48}
EG_T(z,\tau)\to\frac12\cdot\Bigl(\frac{m^{-\xi}-1}{-\xi}\Bigr)\cdot z'\cQ_Hz,
\end{equation}
while Lemma~\ref{lem:9.6} shows that $\Var(G_T(z,\tau))\to0$. In what follows, $F_t$, $f_t$
and $E_t$ denote $F_U(\cdot|X_t)$, $f_U(\cdot|X_t)$ and $E[\cdot|X_t]$, respectively, where
$U$ is the auxiliary error constructed in Condition~R1.

Since
\begin{equation}\label{e:9.49}
\begin{aligned}
G_T(z,\tau)&\equiv\sum_{t=1}^Ta_T\cdot\Bigl(\int_0^{X_t'z/a_T}\Bigl[\frac{\ind(Y_t-X_t'\beta(\tau)\le s)-\ind(Y_t-X_t'\beta(\tau)\le0)}{\sqrt{\tau T}}\Bigr]ds\Bigr)\\
&=\sum_{t=1}^T\Bigl(\int_0^{X_t'z}\Bigl[\frac{\ind(Y_t-X_t'\beta(\tau)\le s/a_T)-\ind(Y_t-X_t'\beta(\tau)\le0)}{\sqrt{\tau T}}\Bigr]ds\Bigr),
\end{aligned}
\end{equation}
we have
\begin{align}
EG_T(z,\tau)&=T\cdot E\Bigl(\int_0^{X_t'z}\frac{F_t[F_t^{-1}(\tau)+s/a_T]-F_t[F_t^{-1}(\tau)]}{\sqrt{\tau T}}\,ds\Bigr)\nonumber\\
&\overset{(1)}{=}T\cdot E\Bigl(\int_0^{X_t'z}\frac{f_t\{F_t^{-1}(\tau)+o(F_u^{-1}(m\tau)-F_u^{-1}(\tau))\}}{a_T\cdot\sqrt{\tau T}}\cdot s\cdot ds\Bigr)\nonumber\\
&\overset{(2)}{\sim}T\cdot E\Bigl(\int_0^{X_t'z}\frac{f_t\{F_t^{-1}(\tau)\}}{a_T\cdot\sqrt{\tau T}}\cdot s\cdot ds\Bigr)\nonumber\\
&=T\cdot E\Bigl(\frac12\cdot(X_t'z)^2\cdot\frac{f_t\{F_t^{-1}(\tau)\}}{a_T\cdot\sqrt{\tau T}}\Bigr)\label{e:9.50}\\
&=E\Bigl(\frac12\cdot(X_t'z)^2\cdot\frac{F_u^{-1}(m\tau)-F_u^{-1}(\tau)}{\tau(f_t\{F_t^{-1}(\tau)\})^{-1}}\Bigr)\nonumber\\
&\overset{(3)}{\sim}E\Bigl(\frac12\cdot(X_t'z)^2\cdot\frac1{H(X)}\cdot\frac{m^{-\xi}-1}{-\xi}\Bigr)\nonumber\\
&\equiv\frac12\cdot\frac{m^{-\xi}-1}{-\xi}\cdot z'\cQ_Hz.\nonumber
\end{align}
Equality (1) is by the definition of $a_T$ and a Taylor expansion. Indeed, since
$\tau T\to\infty$ uniformly over $s$ in any compact subset of $\R$,
\begin{equation}\label{e:9.51}
s/a_T=s\cdot(F_u^{-1}(m\tau)-F_u^{-1}(\tau))/\sqrt{\tau T}=o(F_u^{-1}(m\tau)-F_u^{-1}(\tau)).
\end{equation}
To show equivalence (2), it suffices to prove that, for any sequence
$v_\tau=o(F_u^{-1}(m\tau)-F_u^{-1}(\tau))$ with $m>1$ as $\tau\searrow0$,
\begin{equation}\label{e:9.52}
f_t(F_t^{-1}(\tau)+v_\tau)\sim f_t(F_t^{-1}(\tau))\qquad\text{uniformly in }t.
\end{equation}
This will be shown by using the assumption made in Condition~R3, which is that uniformly
in $t$, $1/f_t(F_t^{-1}(\tau))\sim\partial F_u^{-1}(\tau/K(X_t))/\partial\tau$, where
$\partial F_u^{-1}(\tau)/\partial\tau$ is regularly varying with index $-\xi-1$.

To be clear, let us first show \eqref{e:9.52} for the special case of $f_t=f_u$ and
$F_t^{-1}(\tau)=F_u^{-1}(\tau)$:
\begin{equation}\label{e:9.53}
f_u(F_u^{-1}(\tau)+v_\tau)\sim f_u(F_u^{-1}(\tau)).
\end{equation}
By the regular variation property of
$\partial F_u^{-1}(\tau)/\partial\tau=1/f_u(F_u^{-1}(\tau))$, locally uniformly in $l$
[uniformly in $l$ in any compact subset of $(0,\infty)$],
\begin{equation}\label{e:9.54}
f_u(F_u^{-1}(l\tau))\sim l^{\xi+1}f_u(F_u^{-1}(\tau)).
\end{equation}
That is, locally uniformly in $l$,
\begin{equation}\label{e:9.55}
f_u(F_u^{-1}(\tau)+[F_u^{-1}(l\tau)-F_u^{-1}(\tau)])\sim l^{\xi+1}f_u(F_u^{-1}(\tau)).
\end{equation}
Hence, for any $l_\tau\to1$,
\begin{equation}\label{e:9.56}
f_u(F_u^{-1}(\tau)+[F_u^{-1}(l_\tau\tau)-F_u^{-1}(\tau)])\sim f_u(F_u^{-1}(\tau)).
\end{equation}
Hence, for any sequence $v_\tau=o([F^{-1}(m\tau)-F^{-1}(\tau)])$ with $m>1$ as
$\tau\searrow0$,
\begin{equation}\label{e:9.57}
f_u(F^{-1}(\tau)+v_\tau)\sim f_u(F_u^{-1}(\tau)),
\end{equation}
because for any such $\{v_\tau\}$, in view of Lemma~\ref{lem:9.2}(iii), we can choose a
sequence $\{l_\tau\}$ such that $\{v_\tau\}=\{[F_u^{-1}(l_\tau\tau)-F_u^{-1}(\tau)]\}$ and
$l_\tau\to1$ as $\tau\searrow0$.

Next, let us strengthen the claim \eqref{e:9.53} to \eqref{e:9.52}, completing the proof of
equivalence (2) in \eqref{e:9.50}. Since
\begin{enumerate}
\item[(a)] $1/f_t(F_t^{-1}(\tau))\sim\partial F_u^{-1}(\tau/K(X_t))/\partial\tau
=1/\{K(X_t)f_u[F_u^{-1}(\tau/K(X_t))]\}$ uniformly in $t$ by Condition~R3, and
\item[(b)] $f_u(F_u^{-1}(l\tau/K))\sim(l/K)^{\xi+1}f_u(F_u^{-1}(\tau))\sim(l)^{\xi+1}f_u(F_u^{-1}(\tau/K))$,
locally uniformly in $l$ and uniformly in $K\in\{K(x)\colon x\in\X\}$ [compact by
assumptions on $K(\cdot)$ and $\X$], by \eqref{e:9.54}
\end{enumerate}
we have that locally uniformly in $l$ and uniformly in $t$,
\begin{equation}\label{e:9.58}
f_t(F_t^{-1}(l\tau))\sim l^{\xi+1}f_t(F_t^{-1}(\tau)).
\end{equation}
Repeating the steps \eqref{e:9.55}--\eqref{e:9.57} with $f_t(F_t^{-1}(l\tau))$ in place of
$f_u(F_u^{-1}(l\tau))$, we obtain the required conclusion \eqref{e:9.52}.

The equivalence (3) in \eqref{e:9.50} can be shown as follows. By (a), uniformly in $t$,
\begin{equation}\label{e:9.59}
\frac{F_u^{-1}(m\tau)-F_u^{-1}(\tau)}{\tau(f_t[F_t^{-1}(\tau)])^{-1}}\sim
\frac{F_u^{-1}(m\tau)-F_u^{-1}(\tau)}{\tau(K(X_t)f_u[F_u^{-1}(\tau/K(X_t))])^{-1}}.
\end{equation}
By (b) we have that uniformly in $t$,
\begin{equation}\label{e:9.60}
f_u[F_u^{-1}(\tau/K(X_t))]\sim(1/K(X_t))^{\xi+1}\cdot f_u(F^{-1}(\tau)).
\end{equation}
Putting \eqref{e:9.59} and \eqref{e:9.60} together, we have uniformly in $t$,
\begin{align}
\frac{F_u^{-1}(m\tau)-F_u^{-1}(\tau)}{\tau(f_t[F_t^{-1}(\tau)])^{-1}}
&\sim\frac1{K(X_t)^\xi}\cdot\frac{F_u^{-1}(m\tau)-F_u^{-1}(\tau)}{\tau(f_u[F_u^{-1}(\tau)])^{-1}}\label{e:9.61}\\
&=\frac1{H(X_t)}\cdot\frac{F_u^{-1}(m\tau)-F_u^{-1}(\tau)}{\tau(f_u[F_u^{-1}(\tau)])^{-1}},\label{e:9.62}
\end{align}
where $H(X_t)=X_t'\bc$ for $\xi\ne0$ and $H(X_t)=1$ for $\xi=0$. Finally, by the regular
variation property, \eqref{e:9.54},
\begin{align}
\frac{F_u^{-1}(m\tau)-F_u^{-1}(\tau)}{\tau(f_u[F_u^{-1}(\tau)])^{-1}}
&\equiv\int_1^m\frac{f_u[F_u^{-1}(\tau)]}{f_u[F_u^{-1}(s\tau)]}\,ds\label{e:9.63}\\
&\sim\int_1^ms^{-\xi-1}\,ds\label{e:9.64}\\
&=\frac{m^{-\xi}-1}{-\xi}\qquad(\ln m\text{ if }\xi=0).\label{e:9.65}
\end{align}
Putting \eqref{e:9.61}--\eqref{e:9.65} together gives (3) in \eqref{e:9.50}.

\medskip\noindent\textsc{Step 3.} For $\hZ_T(l)$ defined in \eqref{e:5.3}, notice that
\begin{align*}
(\hZ_T(l_i),\ i=1,\dots,k)&\in\argmin_{z\in\R^{d\times k}}[Q_T(z_1,l_1\tau)+\cdots+Q_T(z_k,l_k\tau)]\\
&=\argmin_{z\in\R^{d\times k}}\Bigl[\sum_{i=1}^kW_T(\tau l_i)'\times z_i+G_T(z_i,\tau l_i)\Bigr]
\end{align*}
for $z=(z_1',\dots,z_k')'$, where the functions $Q_T(\cdot,\cdot)$, $W_T(\cdot)$ and
$G_T(\cdot,\cdot)$ are defined in \eqref{e:9.45}. Since this objective function is a sum of
the objective functions in the preceding steps, it retains the properties of the elements
summed. Therefore, the previous argument applies to conclude that the marginal limit of
this objective function is given by $\sum_{i=1}^kW(l_i)'z_i+G(z_i,l_i)$, where
$(W(l_i),\ i\le k)\equiv N(0,\Sigma)$ with
$EW(l_i)W(l_j)'\equiv EXX'\min(l_i,l_j)/\sqrt{l_il_j}$ and, by calculations that are
identical to those in the preceding section,
$G(z,l_i)\equiv G(z)\equiv\frac12\cdot(\frac{m^{-\xi}-1}{-\xi})\cdot z'\cQ_Hz$. The limit
objective function is minimized at
$(Z_\infty(l_i),\ i\le k)=(\frac{\xi}{m^{-\xi}-1}\cdot\cQ_H^{-1}W(l_i),\ i\le k)$.
Therefore, $(\hZ_T(l_i),\ i\le k)\dto(Z_\infty(l_i),\ i\le k)$.

\medskip\noindent\textsc{Step 4.} It suffices to prove the result for $l=1$. Then
\begin{equation}\label{e:9.66}
\begin{aligned}
\frac{\bar X'(\hat\beta(m\tau)-\hat\beta(\tau))}{\mu_X'(\beta(m\tau)-\beta(\tau))}
&\equiv\frac{\bar X'(\hat\beta(m\tau)-\beta(m\tau))}{\mu_X'(\beta(m\tau)-\beta(\tau))}\\
&\quad-\frac{\bar X'(\hat\beta(\tau)-\beta(\tau))}{\mu_X'(\beta(m\tau)-\beta(\tau))}
+\frac{\bar X'(\beta(m\tau)-\beta(\tau))}{\mu_X'(\beta(m\tau)-\beta(\tau))}\pto1,
\end{aligned}
\end{equation}
since the first two elements on the right-hand side are
$O_p(\tfrac1{\sqrt{\tau T}})=o_p(1)$ by the first part of Theorem~\ref{thm:5.1}.
\hfill$\square$

% ----------------------------------------------------------------------
\subsection{CLT for $W_T(\tau)$ and LLN for $G_T(z,\tau)$}\label{sec:9.6}

\begin{lemma}[CLT and LLN]\label{lem:9.6}
Let $\{Y_j,X_j\}_{-\infty}^t$ be an i.i.d.\ or a stationary $\alpha$-mixing sequence. The
following statements are true for $W_T(\cdot)$ and $G_T(\cdot,\cdot)$, defined in
\eqref{e:9.45}, as $\tau\searrow0$ and $\tau T\to\infty$:
\begin{enumerate}
\item[(i)] Suppose mixing coefficients satisfy $\alpha_j=O(j^{-\phi})$ with $\phi>2$, and
for any $K$ sufficiently close to $0^+$ or $-\infty$, uniformly in $t$ and $s\ge1$, and
some $C>0$,
\begin{equation}\label{e:9.67}
P(U_t\le K,U_{t+s}\le K)\le CP(U_t\le K)^2.
\end{equation}
Then for any finite collection of positive constants $l_1,\dots,l_k$,
\[
\{W_T(\tau l_1)',\dots,W_T(\tau l_k)'\}'\dto(W(l_1)',\dots,W(l_k)')'=N(0,\Sigma)
\]
with $EW(l_i)W(l_j)'\equiv EXX'\min(l_i,l_j)/\sqrt{l_il_j}$.
\item[(ii)] If, in addition, $\alpha_j=O(j^{-\phi})$ with $\phi>\frac1{1-\gamma}$ for
$0<\gamma<1$ and $\tau^{1-2/\gamma}/T\to0$, then
\begin{equation}\label{e:9.68}
\Var(G_T(z,\tau))\to0.
\end{equation}
\end{enumerate}
\end{lemma}

\begin{remark}\label{rem:9.3}
In the i.i.d.\ case the claim (i) simply follows from the Lindeberg--Feller CLT. In the
dependent case condition \eqref{e:9.67} requires that the extremal events should not
cluster, which leads to the same limits as under i.i.d.\ sampling. This condition may
possibly be refined along the lines of \citet{WRL}, who dealt with the nonregression case.
\eqref{e:9.67} is analogous to the no-clustering conditions of Robinson
[\citeyear{Robinson}, A7.4, page 191] used in the context of kernel estimation.
\end{remark}

\begin{proof}[Proof of Lemma~\ref{lem:9.6}]
To show (i), $\{W_T(\tau l_i)',\ i\le k\}'$ suits the CLT of \citet{Robinson}, which
implies the same weak limit as under i.i.d.\ sampling. His conditions A7.1 (with $q=0$),
A7.2 and A7.3 are satisfied automatically. The assumed above mixing condition implies
$\sum_{j=1}^\infty j\alpha_j<\infty$, which implies his condition A3.3. Last, condition
\eqref{e:9.67} immediately implies his condition A7.4.

To show (ii), suppress $\tau$. Then from \eqref{e:9.49},
\[
\Var(G_T(z))=\tau^{-1}\Bigl(\Var(\lambda_1)+2\sum_{k=1}^{T-1}\frac{T-k}{T}\Cov(\lambda_1,\lambda_{1+k})\Bigr),
\]
for
\[
\lambda_t=\int_0^{X_t'z}\bigl[\ind(Y_t-X_t'\beta(\tau)\le s/a_T)-\ind(Y_t-X_t'\beta(\tau)\le0)\bigr]\,ds.
\]
By Condition~R2, $|\lambda_t|\le K_0|\mu_t|$, for
\[
\mu_t=\bigl(\ind(Y_t-X_t'\beta(\tau)\le X_t'z/a_T)-\ind(Y_t-X_t'\beta(\tau)\le0)\bigr)
\]
and some $K_0<\infty$. Hence,
\begin{equation}\label{e:9.69}
\begin{aligned}
\Var(\lambda_1)&=O(E\lambda_1^2)\overset{(1)}{=}O(E\mu_1^2)\overset{(2)}{=}O(E|\mu_1|)\\
&\overset{(3)}{=}O(f_u(F_u^{-1}(\tau))a_T^{-1})=O\bigl(\sqrt{\tau/T}\bigr),
\end{aligned}
\end{equation}
where (1) is by $|\lambda_t|\le K_0|\mu_t|$, (2) is by $|\mu_t|\in\{0,1\}$, and (3) is by
the calculation in \eqref{e:9.50}. Thus, in the i.i.d.\ case $\Var(G_T(z))=o(1)$ follows
from \eqref{e:9.69} and $\tau T\to\infty$. Also, for all $s$ and some positive constants
$K_1,K_2,K_3,K_4$,
\begin{align}
|\Cov(\lambda_1,\lambda_{1+s})|&\le K_1\bigl(\alpha_s^{1-\gamma}[E|\lambda_1|^r]^{1/r}[E|\lambda_1|^p]^{1/p}\bigr)\label{e:9.70}\\
&\le K_2\bigl(\alpha_s^{1-\gamma}[E|\mu_1|^r]^{1/r}[E|\mu_1|^p]^{1/p}\bigr)\label{e:9.71}\\
&\le K_3\bigl(\alpha_s^{1-\gamma}[E|\mu_1|]^{\gamma}\bigr)\label{e:9.72}\\
&\le K_4\Bigl(\alpha_s^{1-\gamma}\Bigl(\frac\tau T\Bigr)^{\gamma/2}\Bigr),\label{e:9.73}
\end{align}
where $1/p+1/r=\gamma\in(0,1)$, $p\ge1$. Here \eqref{e:9.70} follows by Ibragimov's mixing
inequality [e.g., \citet{Davidson}], \eqref{e:9.71} follows by the previous bound
$|\lambda_t|\le K_0|\mu_t|$ and \eqref{e:9.72} follows by $|\mu_t|\in\{0,1\}$, while
\eqref{e:9.73} follows by \eqref{e:9.69}. So $\Var(G_T(z))=o(1)$ by the condition on the
mixing coefficients.
\end{proof}

% ----------------------------------------------------------------------
\subsection{Proof of Theorem~\ref{thm:6.1}}\label{sec:9.7}
Theorem~\ref{thm:6.1} is a direct corollary of Theorem~\ref{thm:5.1} and
Lemma~\ref{lem:9.1}. Proof of claim (i) follows similarly to the proof in \eqref{e:9.66}.
Claim (i) implies claims (ii)--(iv), using the properties (v) and (vi) in
Lemma~\ref{lem:9.1}. Uniformity in $x$ in claim (iv) follows from the linearity of
$\hat\rho_{x,\bar X,1}$ in $x$. Finally, claim (v) follows from Theorem~\ref{thm:5.1}
by the delta method.

% ----------------------------------------------------------------------
\subsection{Tightness of $Z_\infty(k)$}\label{sec:9.8}
This section provides primitive conditions for tightness of $Z_\infty(k)$, which is assumed
in the statement of Theorem~\ref{thm:4.1} and the conditions of uniqueness given in
Remark~\ref{rem:4.4}.

We impose the design condition of \citet{PJ}, who used it for the case $\tau T\to0$ and
show its plausibility on page 233, for example, when $EXX'>0$. Their proof of tightness is
not applicable here, so we have it.

\begin{conditionPJ}\label{PJ}
Let $F_X$ denote the distribution function of $X$. There are a finite integer $I$, a
collection of sets $\{R_1,\dots,R_I\}$ and positive constants $\delta$ and $\eta$ such
that:
\begin{enumerate}
\renewcommand{\labelenumi}{(\alph{enumi})}
\item for each $u\in\{u\colon\|u\|\ge1,\ u_1\ge0\}$, there is $R_{j(u)}$ such that
$x'u>\delta\|u\|$ for all $x\in R_{j(u)}$,
\item $\int_{R_j}dF_X>\eta>0$ for all $j=1,\dots,I$.
\end{enumerate}
\end{conditionPJ}

\begin{lemma}\label{lem:9.7}
If Conditions R1, R2 and PJ hold, then $Z_\infty(k)$ is finite a.s.
\end{lemma}

\begin{proof}
Choose $z^f=(z_1^f,\dots,z_d^f)'\in\R^d$ such that
\begin{equation}\label{e:9.74}
Q_\infty(z^f,k)\equiv-k\mu_X'z^f+\int_E(x'z^f-u)^+\,d\bN(u,x)=O_p(1),
\end{equation}
which is possible, as shown in the proof of Theorem~\ref{thm:4.1}.

Consider a closed ball $B(M)$ with radius $M$ and center $z^f$, and let
$z(k)=z^f+\delta(k)v(k)$, where $v(k)=(v_1(k),\dots,v_d(k))'$ is a direction vector with
unity norm $\|v(k)\|=1$ and $\delta(k)\ge M$. By convexity in $z$,
\begin{equation}\label{e:9.75}
\frac M{\delta(k)}\bigl(Q_\infty(z(k),k)-Q_\infty(z^f,k)\bigr)\ge Q_\infty(z^*(k),k)-Q_\infty(z^f,k),
\end{equation}
where $z^*(k)$ is a point of boundary of $B(M)$ on the line connecting $z(k)$ and $z^f$. We
will show that, for any $K$ and $\varepsilon>0$, there is large $M$ such that
\begin{equation}\label{e:9.76}
P\Bigl(\inf_{v(k)\colon\|v(k)\|=1}Q_\infty(z^*(k),k)>K\Bigr)\ge1-\varepsilon.
\end{equation}
\eqref{e:9.76} and \eqref{e:9.74} imply \eqref{e:9.75} $>C>0$ with probability arbitrarily
close to $1$ for $M$ sufficiently large, meaning that $Z_\infty(k)\in B(M)$ with
probability arbitrarily close to $1$ for $M$ sufficiently large, that is,
$Z_\infty(k)=O_p(1)$.

Thus, it remains to show \eqref{e:9.76}. Since $\mu_X=(1,0,\dots,0)'$,
$\mu_X'z^*(k)=z_1^f+v_1(k)\cdot M$. Hence, it suffices to show that, for any
$\varepsilon>0$ and any large $K>0$,
\begin{equation}\label{e:9.77}
\begin{gathered}
-v_1(k)\cdot k\cdot M+\int_E(x'z^*(k)-u)^+\,d\bN(u,x)\ge K\\
\text{w.pr.\ }\ge1-\varepsilon,\text{ for large enough }M,
\end{gathered}
\end{equation}
and, therefore, we establish \eqref{e:9.76}. Fix $\epsilon_0\in(0,\delta/3)$ and set
$\tilde\delta=\delta/2$. If $v_1(k)\le-\epsilon_0$, then since the integral in
\eqref{e:9.77} is nonnegative, the left-hand side of \eqref{e:9.77} is at least
$\epsilon_0 kM$, which exceeds $K$ for $M$ sufficiently large. Now suppose
$v_1(k)>-\epsilon_0$. If $v_1(k)\ge0$, Condition~PJ yields $R_{j(v)}$ with
$x'v(k)>\delta>\tilde\delta$ for all $x\in R_{j(v)}$; if $-\epsilon_0<v_1(k)<0$,
applying Condition~PJ to $\tilde v=(0,v_{-1}(k)')'/\|v_{-1}(k)\|$ yields $R_{j(\tilde v)}$
on which $x'v(k)=v_1(k)+\|v_{-1}(k)\|x'\tilde v>-\epsilon_0+\sqrt{1-\epsilon_0^2}\,\delta
>\tilde\delta$. Thus, for every $v(k)$ with $v_1(k)>-\epsilon_0$, there is $j(v)\in\{1,\dots,I\}$
such that
\begin{equation}\label{e:9.78}
\begin{aligned}
&\int_E(x'z^*(k)-u)^+\,d\bN(u,x)\\
&\quad\ge\int_{([-\infty,\kappa]\times R_{j(v)})\cap E}(x'z^*(k)-u)^+\,d\bN(u,x)\\
&\quad\ge\bN(([-\infty,\kappa]\times R_{j(v)})\cap E)\times(\tilde\delta M-\kappa-\kappa')^+,
\end{aligned}
\end{equation}
where $\kappa\in\R$ is a constant to be determined later and that does not depend on
$v(k)$ and $\kappa'=\max_{x\in\X}|x'z^f|$.

Note that for any region $\mathbb X$ such that $\int_{\mathbb X}dF_X>\eta>0$ and any
$\kappa_1>0$ and $\varepsilon>0$, there exists $\kappa_2$ (sufficiently large for type 1 and 3
tails, and negative and sufficiently close to $0^-$ for type 2 tails) such that
\begin{equation}\label{e:9.79}
\bN(([-\infty,\kappa_2]\times\mathbb X)\cap E)>\kappa_1\qquad\text{w.pr.\ }\ge1-\varepsilon.
\end{equation}
Hence, by \eqref{e:9.79} we can select $\kappa$ (sufficiently large for type 1 and 3 tails,
and in $(-1,0)$ sufficiently close to $0^-$ for type 2 tails) so that
\begin{equation}\label{e:9.80}
\begin{gathered}
\bN(([-\infty,\kappa]\times R_j)\cap E)>\frac{(k+1)}{\tilde\delta}\\
\text{for all }j\in\{1,\dots,I\}\ \text{w.pr.\ }\ge1-\varepsilon,
\end{gathered}
\end{equation}
so that uniformly over $v(k)$ with $v_1(k)>-\epsilon_0$, w.pr.\ $\ge1-\varepsilon$,
\begin{equation}\label{e:9.81}
\begin{aligned}
&-v_1(k)\cdot k\cdot M+\int_E(x'z^*(k)-u)^+\,d\bN(u,x)\\
&\qquad\ge-k\cdot M+(k+1)\frac{(\tilde\delta M-\kappa-\kappa')^+}{\tilde\delta}.
\end{aligned}
\end{equation}
Now set $M$ sufficiently large to obtain \eqref{e:9.77} uniformly over $\|v(k)\|=1$.
\end{proof}

\section*{Acknowledgments}
This paper is based on a chapter of my Stanford Ph.D.\ dissertation completed in August
2000, and I would like to thank my dissertation committee: Takeshi Amemiya, Pat Bajari and
Tom Macurdy. I would also like to thank Roger Koenker, Stephen Portnoy, Ivan
Fernandez-Val, Jerry Hausman, Bo Honor\'e, Jana Jure\v{c}kov\'a and Keith Knight. I am also
very grateful to two anonymous referees, an Associate Editor and Co-Editors of the journal
for providing feedback of the highest quality.

\bigskip
\noindent
Department of Economics\\
Massachusetts Institute of Technology\\
50 Memorial Drive\\
Cambridge, Massachusetts 02142\\
USA\\
e-mail: \texttt{vchern@MIT.edu}

\end{document}